\documentclass[12pt]{article}

\usepackage{latexsym,amssymb,amsmath}
\usepackage[dvips]{graphicx}
\usepackage{bm}
\newcommand{\Z}{\mathbb Z}
\newcommand{\N}{\mathbb N}
\newcommand{\R}{\mathbb R}
\newcommand{\Q}{\mathbb Q}
\newcommand{\HH}{\mathbb H}

\newcommand{\sma}{\left(\begin{array}}
\newcommand{\fma}{\end{array}\right)}

\newtheorem{lem}{Lemma}[section]
\newtheorem{defn}[lem]{Definition}

\newtheorem{co}[lem]{Corollary}
\newtheorem{thm}[lem]{Theorem}
\newtheorem{prop}[lem]{Proposition}

\newenvironment{proof}{\textbf{Proof.}}{\newline\hspace*{\fill}{$\Box$}\\}

\begin{document}
\title{Groups acting on hyperbolic  spaces with a locally finite orbit}
\author{J.\,O.\,Button\\
Selwyn College\\
University of Cambridge\\
Cambridge CB3 9DQ\\
U.K.\\
\texttt{j.o.button@cam.ac.uk}}
\date{}
\maketitle
\begin{abstract}
A group with a geometric action on some hyperbolic space is
necessarily word hyperbolic, but on the other hand
every countable group acts (metrically) properly by isometries
on a locally finite hyperbolic graph.
In this paper we consider what happens when a group acts isometrically
on a restricted class of hyperbolic spaces, for instance quasitrees.
We obtain strong conclusions on the group structure if the action has
a locally finite orbit, especially if the group is finitely generated.

We also look at group actions on finite products of quasitrees, where
our actions may be by automorphisms or by isometries, including the
Leary - Minasyan group.
\end{abstract}

\section{Introduction}

In this paper we consider groups (these are abstract but not
necessarily finitely generated and they may even be uncountable) acting
by isometries on geodesic metric spaces $X$, where $X$ will generally be
a hyperbolic space or a finite product of hyperbolic spaces. If such an
action is proper and cocompact then $X$ is a proper metric space (all
closed balls are compact) and $G$ will be finitely generated. Thus we
can put the word metric on $G$ and then $G$ and $X$ are quasi-isometric
by the Svarc - Milnor Theorem. But this theorem further gives us a
quasi-isometry $q:G\rightarrow X$ that is obtained from the group structure
of $G$ and not just its structure as a metric space: namely $q$ is the
orbit map $g\mapsto g(x_0)$ for $x_0\in X$ an arbitrary base point.

Moreover if $G$ is finitely generated and the space $X$ is proper then
the converse holds too: if the orbit map is a quasi-isometry then the
action is proper and cocompact. We might thus wonder
what happens if $X$ is not
proper and also what if the orbit map is not a quasi-isometry but still
well behaved geometrically. In the former case we do obtain the converse
for an arbitrary geodesic metric space $X$ if we replace proper and
cocompact by metrically proper and cobounded. Note though that in
general the orbit map is quasi-surjective if and only if the action
is cobounded. This raises the question: what can we say if a finitely
generated group $G$ acts isometrically on a geodesic metric space $X$
and the orbit map is just a quasi-isometric embedding, rather than a full
quasi-isometry? This does imply that the action is metrically proper
but is in fact strictly stronger. For instance every countable group
acts (metrically) properly on a locally finite hyperbolic graph, by using
the hyperbolic cone construction, but if $G$ has a quasi-isometric embedding
in any hyperbolic space then it is also hyperbolic.

This raises the question of whether we can impose further conditions on
the group $G$ or the hyperbolic space $X$, or the action of $G$ on $X$,
in order to obtain some meaningful conclusion about the structure of
$G$. We can also ask: if there is an action of $G$ on some
space $X$ of a given type, say (metrically) proper, acylindrical or
unbounded, then does $G$ have the same type of action on some space
$Y$ in a more restrictive class? For instance we could take $Y$ to be a
simplicial tree and $X$ to be a quasitree (any geodesic metric space
that is quasi-isometric to a simplicial tree). Here the ``same type''
can be interpreted in two ways: either the existence of some action of $G$
on $Y$ with the same properties, or that this action is directly
related to the given action of $G$ on $X$. For this latter version, the
notion of coarse $G$-equivariant quasi-isometry from one action of $G$
to another is used in \cite{abos}, giving rise to an equivalence class
of actions of a group on different hyperbolic spaces. This has the
advantage that many properties holding for the first action automatically
carry over to the other, so if such a quasi-isometry exists then
we would have a yes answer for our first interpretation of the ``same type''
too.

However \cite{abos} looks at the above notion of equivalence amongst
cobounded actions only. But being cobounded is a property preserved
by coarse $G$-equivariant quasi-isometries, so if we ask: given
an action of $G$ on a hyperbolic space $X$ which is not cobounded, when
is there an action of $G$ on some other hyperbolic space
$Y$ that is cobounded and such that there is a coarse $G$-equivariant
quasi-isometry between them, the answer is never.

Therefore we could replace ``quasi-isometry'' by quasi-isometric
embedding (as was done in \cite{mn3}),
though a coarse $G$-equivariant quasi-isometric embedding
from $X$ to $Y$ also preserves coboundedness. However we will be
going the other way, so we will be looking for a hyperbolic space $Y$
and an action of $G$ on $Y$ such that there is a coarse $G$-equivariant
quasi-isometric embedding from $Y$ to $X$, whereupon a non cobounded
action of $G$ on $X$ can become a cobounded action on $Y$ (and this
map will necessarily not be a quasi-isometry). In many of our applications
$X$ will be some hyperbolic space but $Y$ will be a connected graph
equipped with the path metric. One way that this could be done from any
geodesic metric space $X$ is to use the Rips graph (or Rips complex)
whose vertex set consists of the points of $X$ and an edge is
placed between two points of distance at most 1 away. Then the action of $G$
on $X$ naturally carries over to an action by automorphisms on the Rips graph,
which will be connected and so this is an action by isometries too. Thus
we have a natural $G$-equivariant map from $X$ to the Rips graph which will
be a quasi-isometry, but this graph could be much bigger than $X$.

Here we will adopt a reverse approach and obtain a graph $Y$ from $X$ which
takes into account the action of $G$ but in the smallest way possible. We
call this the Rips orbit graph and it is constructed by choosing a point
$x_0\in X$ and $r>0$ (the result could well depend on $x_0$ and $r$) and
forming the graph with vertex set the orbit of $x_0$ under $G$ and an edge
between two vertices if the distance between them is at most $r$. Then $G$
acts on this graph by automorphisms and thus by isometries if we can find
a value of $r$ for which this graph is connected, which always exists for
finitely generated groups. Moreover we have a $G$-equivariant map from the
vertices of the Rips orbit graph to $X$ and thus a coarse $G$-equivariant
map from the Rips orbit graph itself to $X$. However this need not be a
quasi-isometric embedding, though various properties of the original action
are preserved.

We cover background material in Sections 2 and 3. We begin in Section 2
by looking at the definitions of the different types of isometric action
on a metric space $X$ that we will need, and the maps between these actions.
This includes the notion in Definition 2.1 of an action having a locally
finite orbit, which will be the essential finiteness property of the action
that we will generally require in our results (sometimes along with the
finiteness property of being finitely generated for our groups). We then
define the Rips orbit graph in Definition 2.2 and look at which properties
of a group action are preserved under it.

In that section it is only assumed that $X$ is a metric space (or occasionally
a geodesic metric space). However our focus in this paper is when $X$ is
hyperbolic (meaning that it is geodesic and satisfies any of the equivalent
definitions of $\delta$ - hyperbolicity). In Section 3 we review definitions
and results concerning isometric actions on hyperbolic spaces, with an
emphasis on arbitrary actions and arbitrary hyperbolic spaces. After first
considering individual elements, we review the Gromov classification of
all actions into 5 types: bounded, parabolic, lineal, quasi-parabolic
and general. This will be of importance to us later as some arguments will
require working
through these different types. We also give various examples, in
particular when $X$ is a simplicial tree where we note that there is no
parabolic action of a finitely generated group on $X$ but there can be if
$G$ is not finitely generated. Indeed here we can have a proper action on
a bounded valence tree which is parabolic.

We come to our first results in Section 4 where we look at lineal actions
(that is where the action has two limit points) of an arbitrary (even
uncountable) group $G$ on any hyperbolic space. We might hope that this
implies the existence of a homomorphism from $G$ (or an index 2 subgroup) to
$\R$ but in general we only obtain a homogeneous quasi-morphism. However
Theorem 4.2 states that if the action of $G$ has a locally finite orbit
then we do have such a homomorphism, indeed to $\Z$ (or to $C_2*C_2)$ with
a kernel that acts with bounded orbits. This proof utilises the Rips orbit
graph, despite the fact that it need not be connected in general if $G$ is
not finitely generated. This gives us in Corollary 4.3 our first positive
result on the problem mentioned above of finding an action of $G$ on a
simpler space $Y$ of the same type as the original action on $X$ in Theorem 4.2.
It states that we can take $Y$ to be the simplicial line, where $G$ acts by
automorphisms on $Y$, and we have a coarse $G$-equivariant quasi-isometric
embedding of $Y$ in $X$.

In Sections 5 and 6 we examine isometric actions on quasitrees. Here a
quasitree is simply a geodesic metric space that is quasi-isometric to
a simplicial tree, with no further conditions imposed. Section 5 begins
with several examples of how statements known to hold for actions of
groups on a tree do not hold on replacing tree with quasitree. For
instance if $G=\langle x,y\rangle$ acts on a tree with $x,y$ and $xy$ all
elliptic then the action has bounded orbits (and indeed a global fixed
point) whereas it can be unbounded on a quasitree. However we then follow
our theme of assuming our finiteness condition on $G$ holds, that is $G$
is finitely generated. Then Corollary 5.3 tells us that for any isometric
action of $G$ on an arbitrary quasitree $Q$, we have a coarse $G$-equivariant
quasi-isometric embedding from the Rips orbit graph $\Gamma_r(x)$ (for large
enough $r$ and arbitrary $x\in Q$) to $Q$. If we further assume our other
finiteness condition, that the action on $Q$ has a locally finite orbit
$G\cdot x_0$ then the graph $\Gamma_r(x_0)$ is of bounded valence and the
action of $G$ on this graph is cobounded. We then extend this in Theorem
5.4 to show that under these two finiteness conditions, the graph is
a quasitree and indeed is quasi-isometric to a simplicial tree which
itself has bounded valence.

This allows a host of corollaries: in particular Corollary 5.6 states that
if a finitely generated group $G$ acts metrically properly on any quasitree
then $G$ is virtually free. Also Corollary 5.8 says that if $G$ acts by
automorphisms on a quasitree which is a locally finite graph (equipped with
the path metric) and this action has finite stabilisers (in particular if
it is a free action) then $G$ is also virtually free. We finish Section 5
with a result that applies to any isometric action of a finitely generated
group $G$ on an arbitrary quasitree: namely the action is never parabolic.
Thus any isometry of a quasitree must be either elliptic or loxodromic
(hence recovering a result due to Manning in \cite{man}) and any isometric
action of $G$ in which all elements are elliptic has bounded orbits
(which complements results of this nature in the very recent preprint
\cite{dmht}).

Previous strong results on converting isometric actions on quasitrees
to trees were in \cite{msw}. This uses the idea of a quasi-action but
in the case above where $G$ acts on the quasitree $Q$ which is quasi-isometric
to the tree $T$ say, we obtain a quasi-action on $T$. Theorem 1 of \cite{msw}
then states that if the action of $G$ (here an arbitrary group) on $Q$ is
cobounded and $T$ is a bounded valence, bushy tree then there is some
other tree $T'$ with the same properties as $T$ and on which $G$ acts by
isometries, along with a coarse $G$-equivariant quasi-isometry from $T$ to
$Q$. This has powerful consequences provided we can satisfy the cobounded
condition for the action on $Q$ and the bounded valence, bushy conditions
on $T$. Our Theorem 6.3 states that if we have our usual two finiteness
conditions of $G$ being finitely generated and the presence of a locally
finite orbit then no further hypotheses are required to conclude that
$G$ acts coboundedly by automorphisms on some bounded valence tree $T'$ and
we obtain a coarse $G$-equivariant quasi-isometric embedding from $T'$ to
$Q$. The proof of this theorem requires our results from Sections 4 and 5
as well as the application of \cite{msw} Theorem 1. In the case where the
quasitree $Q$ is a locally finite graph and the action of $G$ on $Q$ is
unbounded, this will also apply to $G$ acting on $T'$, so Corollary 6.4
concludes that $G$ splits non trivially by standard Bass - Serre theory.

In Sections 7 and 8 we take a brief detour to acylindrical actions and
also the case where the hyperbolic space $X$ is also a bounded valence
(as opposed to locally finite) hyperbolic graph with its path metric.
An acylindrically hyperbolic group has a non elementary acylindrical
action on some hyperbolic space $X$ and many examples of such groups are
known. However, if we ask what groups $G$ have such an action on a
proper hyperbolic space then Theorem 7.1 says that this action must be
proper, so for instance $G$ cannot contain $\Z\times\Z$. Moreover
(Corollary 7.4) if further $G$ is finitely generated and $X$ is a quasitree
then $G$ is virtually free. This contrasts vastly with \cite{balas} which
showed that any acylindrical hyperbolic group acts acylindrically on
some graph which is a quasitree (but not locally finite).

Proposition 8.1 can be seen as a Tits alternative for proper actions of an
arbitrary group on any bounded valence hyperbolic graph: the group must
be virtually cyclic, locally finite or contain $F_2$, as opposed to
locally finite hyperbolic graphs where we have mentioned that all countable
groups can occur.

In the last two sections actions of groups on finite products of hyperbolic
spaces, specifically trees and quasitrees, are investigated. In the case
of a product of trees, any group with a proper action by automorphisms
either contains $F_2$ or is virtually a locally finite by $\Z^n$ group
by Theorem 10.1, and this exact statement carries over to a product
of quasitrees if they are locally finite graphs. However we can also
replace the assumption of an action by automorphisms with that of an
action by isometries which can be more general. In \cite{bebf} a finitely
generated group $G$ is said to have property (QT) if there is an action
of $G$ on a finite product of graphs which are quasitrees (not necessarily
locally finite) such that the orbit map is a quasi-isometric embedding.
To finish we investigate the Leary - Minasyan group in \cite{lrmn} which
acts geometrically on $T\times \R^2$ where $T$ is a regular tree of degree
10. Thus this group has property (QT). Nevertheless we show in Theorem 10.4
that it does not. This dichotomy comes about because the definition in
\cite{bebf} uses the $l_1$ product metric, whereas the $l_2$ product
metric is implicitly being used in the above geometric action. This theorem
is established by looking at the possible stable translation lengths of
$\Z^2$ when acting on an arbitrary hyperbolic space.

The author would like to acknowledge helpful conversations with Jason
Manning (in person, just before the first lockdown) and Alice Kerr
(not in person, during the $n$th lockdown).

\section{Types of group action}

In this paper we will be considering
an arbitrary group $G$ acting on an arbitrary metric space
$(X,d)$ by isometries. We say that $G$ acts {\bf freely} on $X$
if all non-identity elements act without fixed points. However when
it comes to the notion of a proper action, various definitions are
in use. Here we will say that $G$ acts ({\bf topologically})
{\bf properly} if for every compact set $C$ of $X$, the set of
elements $\{g\in G\,|\,g(C)\cap C\neq\emptyset\}$ is finite. This
definition can be applied to any group acting by homeomorphisms on a
topological space, hence our name.

In the case where $G$ is acting by isometries on an arbitrary
metric space $X$ (but not in general if the action is by homeomorphisms),
this definition is equivalent to various other definitions in the
literature, including:\\
(i) The map $(g,x)\mapsto (g(x),x)$ from $G\times X$ to $X\times X$
is a proper map (inverse images of compact sets are compact), where
$G$ is given the discrete topology.\\
(ii) For all compact sets $A,B$ of $X$, the set $\{g\in G\,|\,g(A)\cap
B\neq \emptyset\}$ is finite.\\
(iii) For all $x,y\in X$ there exists an open neighbourhood $U$ of $x$
and $V$ of $y$ such that the set $\{g\in G\,|\, V\cap g(U)\neq \emptyset \}$
is finite.\\
(iv) For all $x\in X$ there exists an open neighbourhood $U$ of $X$
such that the set $\{g\in G\,|\,U\cap g(U)\neq\emptyset \}$ is finite.\\
(v) We do not have a point $x\in X$ and a sequence of distinct elements
$g_n\in G$ with $g_n(x)$ tending to a point in $X$.\\
(vi) We do not have a point $x\in X$ and a sequence of distinct elements
$g_n\in G$ with $g_n(x)$ tending to $x$.

There is another characterisation of a topologically
proper action in a special case
that will be of interest to us throughout this paper. Suppose that
$X$ is a connected combinatorial graph (but not necessarily locally
finite) which is given the path metric
where all edges have length 1. If a group $G$ acts by graph automorphisms 
then it also acts by isometries and in the former case we can see that
acting topologically properly is equivalent
to all vertex stabilisers being finite (eg by using (iv) above).
However this is not quite true
in the latter case: by considering vertices of degree 2, we note that
an isometry of a graph $\Gamma$ equipped with its path metric is indeed
a graph automophism if $\Gamma$ is not the simplicial line or a circle
of $n\geq 1$ edges. But the action of $\Z^2$ generated by
$x\mapsto x+1$ and $x\mapsto x+\sqrt{2}$ on the simplicial line
(or of $\Z$ on a circle via an irrational rotation) is
isometric and free but not by graph automorphisms and not proper. 
In this paper we will sometimes need to keep the distinction between an action
on a graph by graph automorphisms and by isometries (this is especially
true when we deal with products of graphs).

Going back now to the case where $G$ acts on an arbitrary metric space
$X$ by isometries, we also have the notion of a
{\bf metrically proper} action, which
is the same as a topologically proper action except that  ``compact'' is
replaced by ``closed and bounded'', or equivalently for some (or any)
$x_0\in X$ and any $R\geq 0$, the set $\{g\in G\,|\,d(g(x_0),x_0)\leq R\}$ is
finite. If $X$ is a {\bf proper} metric space, that is all closed balls
are compact (for instance any locally finite connected graph), then these two
notions are clearly equivalent but in general acting metrically properly
is somewhat stronger. Similarly we say that an action is {\bf cocompact} 
if there is a compact set whose translates under $G$ cover all of $X$ and
{\bf cobounded} if there is a bounded set. Once again, these are equivalent
concepts in a proper metric space.

If a group $G$ has an isometric action on a metric space
$(X,d_X)$ and $G$ also acts on a set $Y$ then we
say that a function $F:Y\rightarrow X$ is {\bf coarse $G$-equivariant}
if there is $M\geq 0$ such that for all $y\in Y$ and all $g\in G$ we
have $d_X(F(g(y)),g(F(y)))\leq M$. For instance if $Y=G$ then (for an
arbitrary basepoint $x_0\in X$) the orbit map $g\mapsto g(x_0)$ is
$G$-equivariant,
so if $Y$ is a Cayley graph $Cay(G)$ of $G$ then we can extend the
orbit map to $Cay(G)$ by sending any point to a nearest vertex $g$ and then to
$g(x_0)$. This is a coarse $G$-equivariant map with $M=1$. 

In our setting $Y$ will also be a metric space
and $G$ will act by isometries on $Y$ too, so
we are particularly interested
in the case where $F$ is a quasi-isometry or more generally a quasi-isometric
embedding. For these two types of function,
it will be helpful to consider whether a given property that holds
for the action of $G$ on $Y$ carries over to the action of $G$
on $X$. In the case when $F$ is a quasi-isometric embedding but not
necessarily a quasi-isometry, we can also ask conversely, and more
importantly for our purposes, whether a
property that holds for the action on $X$ continues to hold for the
action on $Y$. It is easy to check that 
acting metrically properly is preserved both ways, but acting topologically
properly need not be preserved in either direction (which is one reason why
metrically proper actions are more appropriate in this paper) and nor is
acting freely or faithfully. As for
cobounded actions, this is preserved under a quasi-isometric embedding
when passing from $X$ to $Y$ (so in both directions for a quasi-isometry)
but not in general from $Y$ to $X$. Indeed a main aim for us on being given
an isometric but not cobounded action of a group $G$ on a space $X$
will be to find an isometric action of $G$ on some space $Y$ with a coarse
$G$-equivariant quasi-isometric embedding from $Y$ to $X$.

In general we have the
Svarc - Milnor lemma, which we use in the following two versions.
Suppose that $G$ is a group acting by isometries on a geodesic metric space
$X$. If the action is topologically proper and cocompact, or if it is    
metrically proper and cobounded then $G$ is finitely generated and for
any $x_0\in X$ the orbit map $g\mapsto g(x_0)$ is a $G$-equivariant
quasi-isometry from $G$ (with the word metric obtained from any finite
generating set) to $X$. This then gives us a coarse $G$-equivariant
quasi-isometry from the Cayley graph of $G$ to $X$ (which are both
geodesic metric spaces). However if the action is proper
but we do not have the cobounded/cocompactness condition, it is certainly
not the case in general that $G$ is finitely generated. Even if it is, the
orbit map need not be a quasi-isometric embedding. However when $G$ is
finitely generated we do have the converse: if the orbit map from $G$
(equipped with the word metric from a finite generating set) to $X$ is a
quasi-isometric embedding then $G$ will act metrically properly. We also
note here that the orbit map being a quasi-isometric embedding is
equivalent to the existence of any coarse $G$-equivariant quasi-isometric
embedding from $G$ to $X$ and the same is true for quasi-isometries, because
the existence
of any coarse $G$-equivariant map from $G$ to $X$ that is quasi-surjective
is equivalent to the orbit map being quasi-surjective,
which is then equivalent to the action of $G$ on $X$ being cobounded.

We further say that $G$ acts on $X$ with {\bf (un)bounded orbits} if there
exists (or not)
a point in $X$ whose orbit under $G$ is bounded, whereupon all orbits
are then bounded (or not). This is also preserved both ways under a
coarse $G$-equivariant quasi-isometric embedding.
Also $G$ acts on $X$ with {\bf finite orbits} if
all points
in $X$ have a finite orbit under $G$ (here it is not enough to check just one
orbit and this need not be preserved in either direction under coarse
$G$-equivariant quasi-isometric embeddings).
We now define a related notion.
\begin{defn} If a group $G$ acts on a metric space $X$ by isometries then
we say $G$ has a {\bf locally finite orbit} if there is $x_0\in X$
such that for all $R\geq 0$ the set 
$\{g(x_0)\,|\,g\in G\mbox{ and }d(g(x_0),x_0)\leq R\}$
is finite.
\end{defn}
Note that it does not matter how big the stabilisers are as long as there
are only finitely many images of $x_0$ in any ball around $x_0$. Here this
property need not hold for other $x\in X$ in place of $x_0$ (and it need
not be preserved in either direction under a coarse $G$-equivariant
quasi-isometric embedding, just as for having a finite orbit).
However there are
two cases of interest to us where this definition is satisfied for any
$x_0\in X$:\\
\hfill\\
(1) When $G$ acts metrically properly on $X$. Indeed acting metrically
properly on $X$ is equivalent to having a locally finite orbit $G\cdot x_0$
where the stabiliser of $x_0$ is finite (it is also equivalent to having a
locally finite orbit and acting topologically properly).
\hfill\\
(2) When $G$ acts by graph automorphisms on a locally finite connected
graph $\Gamma$ equipped with the path metric.

For a group $G$ acting by isometries on a metric space, we would like to
obtain from this an action of $G$ on a graph which preserves many of
the essential features of the original action (as discussed above in the
case of coarse $G$-equivariant quasi-isometries/quasi-isometric embeddings)
but which is more tractable to study.
Given a metric space $(X,d_X)$ and $r\geq 0$, we have the {\bf Rips graph}
$\Gamma_r(X)$ which is a simplicial graph with vertex set $X$ and an edge
between $x,y\in X$ if and only if $0<d_X(x,y)\leq r$. We then have
that if $G$ acts on $X$ by isometries then it also acts by automorphisms
on $\Gamma_r(X)$ with the restriction to the vertices giving us the
original combinatorial
action. This is a useful construction but it might not give us
exactly what we want. For instance, we would like to put the path metric
$d_\Gamma$ on $\Gamma_r(X)$, whereupon $G$ will act on $(\Gamma_r(X),d_\Gamma)$
by isometries, but $\Gamma_r(X)$ might not be a connected graph. Even if so,
$\Gamma_r(X)$ might not be quasi-isometric to $X$, for instance $X$ could
be bounded and $\Gamma_r(X)$ unbounded, although if $X$ is
a geodesic metric space then the map sending the points of $X$ to the
vertices of $\Gamma_r(X)$ is a $G$-equivariant quasi-isometry.
Finally it might not be necessary to make this
move: for instance if $X$ were already a connected combinatorial graph 
with the path metric which was locally finite (but not a single point)
then we would be replacing it with a graph where every vertex had uncountable
valence.
  
Our approach here is that, because we want to obtain information about
the group $G$ from the space $X$ on which it acts, we adapt the Rips
construction so that we only look at one orbit of the action of $G$.
This has the advantage of ``shrinking'' $X$ to make things more
manageable,
although it could be that we get a different result depending on which
orbit we use.
\begin{defn} \label{rog}
If $(X,d_X)$ is a metric space with $x_0\in X$ and $G$ is a group acting by 
isometries on $X$ then for any $r\geq 0$ the {\bf Rips orbit graph} 
$\Gamma_r(x_0)$ is the combinatorial graph whose vertices are the points
of the orbit $G\cdot x_0$ and there is an edge between
distinct vertices $v,w\in G\cdot x_0$ if and only if $d_X(v,w)\leq r$. 
\end{defn}
We then have that $G$ acts on the vertices of $\Gamma_r(x_0)$ in the same
way, at least combinatorially, as it did on the orbit of $x_0$ in $X$.
This clearly extends to an action of $G$ on $\Gamma_r(x_0)$
by graph automorphisms which is now vertex transitive.
Although we have not assumed in the definition that $G$ has a
locally finite orbit on $X$, the key point is that if $G\cdot x_0$
is locally finite then $\Gamma_r(x_0)$ is a locally finite graph.
But then as the action of $G$ is vertex transitive, we get in this
case that $\Gamma_r(x_0)$ actually has bounded valence.

Of course the next step should be to put the path metric $d_\Gamma$
on $\Gamma_r(x_0)$, whereupon this action of $G$ would be by isometries
and cobounded.
However it is possible
that $\Gamma_r(x_0)$ fails to be connected for any $r$. 
When $G$ is finitely generated though, that does not happen.
This is because on taking any finite generating set $S$
for $G$ and setting $r=r(S)$ to be the maximum of 
$d_X(s(x_0),x_0)$ for $s\in S$,
we have that $\Gamma_r(x_0)$ is connected by induction on word length:
since \[d_X(s_1^{\epsilon_1}\ldots s_{n+1}^{\epsilon_{n+1}}(x_0), 
s_1^{\epsilon_1}\ldots s_n^{\epsilon_n}(x_0))=d_X(s_{n+1}(x_0),x_0)\leq r\]
for any $s_i\in S$ and $\epsilon_i\in\{\pm 1\}$,
if $s_1^{\epsilon_1}\ldots s_n^{\epsilon_n}(x_0)$ is joined to $x_0$
by a path in $\Gamma_r(x_0)$ then so is
$s_1^{\epsilon_1}\ldots s_{n+1}^{\epsilon_{n+1}}(x_0)$. 

We thus have
$d_X(v,w)\leq rd_\Gamma(v,w)$ for any vertices $v,w$. We might further
hope that there is $c>0$ such that $cd_\Gamma(v,w)\leq d_X(v,w)$, which
would mean that the map sending the vertices of $\Gamma_r(x_0)$ back
to $X$ is a $G$-equivariant quasi-isometric embedding. If so then the map from
$\Gamma_r(x_0)$ to $X$ which is equal to this on the vertices and sending edges
to a nearest vertex and then to $X$ would be a coarse $G$-equivariant
quasi-isometric embedding and hence, as mentioned above,
many features of the two actions
of $G$ on $\Gamma_r(x_0)$ and on $X$ would be the same.
In fact even for finitely generated groups there might not be any
quasi-isometric embedding of the Rips orbit graph in $X$
as the example in Proposition \ref{pres} indicates.
This looks like a serious disadvantage of our construction, but we will
see in the next sections that there are cases where we can show that
we do have a coarse $G$-equivariant quasi-isometric embedding
from the Rips orbit graph to our space $X$, even when $X$
is a hyperbolic space.

However we now show that there are some properties of actions
of a group which will be preserved in general
when moving from $X$ to the Rips orbit graph. If we have
a locally finite orbit then these properties will be preserved the other
way too.

\begin{prop} \label{pres}
Let $G$ be a group acting on the metric space $(X,d_X)$ 
by isometries and let $r>0$
be such that the Rips orbit graph $\Gamma_r(x_0)$ is
connected with path metric $d_\Gamma$ (such an $r$ is guaranteed if
$G$ is finitely generated). If the action of $G$ on $X$
has unbounded orbits/ is metrically proper then
so is the corresponding action of $G$ on $\Gamma_r(x_0)$.
This becomes an if and only if provided that
$G\cdot x_0$ is a locally finite orbit on $X$, but not in general
even if $\Gamma_r(x_0)$ has bounded valence.
\end{prop}
\begin{proof}
As the action on $\Gamma_r(x_0)$
is vertex
transitive, this orbit is bounded if and only if $\Gamma_r(x_0)$ 
is a bounded graph, in which case so is the original orbit of $x_0$
on $X$, and thus all of orbits of $G$ on $X$.

Now suppose the action of $G$ on $X$ has the locally
finite orbit $G\cdot x_0$. Then this orbit is bounded in $X$ if and
only if it is finite, so if and only if $\Gamma_r(x_0)$ is a finite
graph. Moreover in this case
$G$ acts metrically properly on $X$ if and
only if the stabiliser of $x_0$ in $X$ is finite.
But this is the same as the stabiliser of
$x_0$ in $\Gamma_r(x_0)$, where the orbit of $x_0$ is still locally
finite by construction. 

However, without the locally finite orbit condition, we could start by
taking the Cayley graph of an infinite group $G$ with respect to a
finite generating set, but then form $X$ by adding a vertex $v_0$ which 
is joined to every other vertex. We then let $G$ act on $X$ by globally
fixing $v_0$ but acting as before on all other vertices.
Then $\Gamma_1(v)$ is the original
Cayley graph if $v\neq v_0$, and so $G$ acts metrically properly
and with unbounded orbits on it, but $X$ is a bounded metric space.
\end{proof}

Note that if a finitely generated group $G$ acts
metrically properly on some metric space $X$ then $G$ acts metrically
properly and coboundedly on the geodesic metric space $\Gamma_r(x_0)$
(for arbitrary $x_0\in X$ and $r$ such that the graph is connected).
Thus by Svarc - Milnor the orbit map from $G$ with its word metric
to $\Gamma_r(x_0)$ is a $G$-equivariant quasi-isometry. Thus if the
natural map sending the vertices of $\Gamma_r(x_0)$ back to $X$ is
a quasi-isometric embedding then by composition we have that the orbit
map from $G$ to $X$ is a quasi-isometric embedding. (Conversely if the
orbit map is a quasi-isometric embedding then we have mentioned that
$G$ acts metrically
properly on $X$, but the former is generally stronger. For instance a
hyperbolic group $G$ acts metrically properly on its Cayley graph $C$, thus
all subgroups do too. But if we restrict this action to a finitely generated
subgroup $H$ then the orbit map from $H$ to $C$ being a quasi-isometric
embedding would imply that $H$ is itself hyperbolic.)
  
\section{Group actions on hyperbolic spaces}

Suppose again for now that $G$ is any group acting by isometries
on any metric space $(X,d)$. We mentioned in the previous section
some definitions and properties for the action of $G$ on $X$. 
Let us now consider how cyclic subgroups might act.

Given any $g\in G$, it is standard to say that
the element $g$ is {\bf elliptic} under the given action if
the subgroup $\langle g\rangle$ has bounded orbits, which is always 
the case if $g$ has
finite order. It is also standard to say that 
the element $g$ is {\bf loxodromic} under the given action (though
this can also be called hyperbolic) if
the subgroup $\langle g\rangle$ embeds quasi-isometrically in $X$
under the (or rather an) orbit map.
This condition can be expressed in various equivalent ways, but our
definition here will be:
for a point $x_0\in X$, let the {\bf stable translation length} be
\[\tau(g):=\mbox{lim}_{n\rightarrow\infty} \frac{d(g^n(x_0),x_0)}{n}.\]
This limit always exists and is independent of $x_0$, so we say that
$g$ is loxodromic if $\tau(g)>0$. Indeed
$\tau(g)$ is always equal to the infimum of this sequence,
so that $g$ being loxodromic is equivalent to saying that
we have $c>0$ such that for all $x\in X$ and $n\in\N$ (or indeed $n\in\Z$)
we have $d(g^n(x),x)\geq nc$.

This leaves elements $g\in G$ which act with unbounded orbits but
where $\tau(g)=0$. At this level of generality, it
does not seem standard to call these {\bf parabolic} elements though we
will do so here. We also note
that an element $g\in G$ is elliptic/loxodromic/parabolic if
and only if $g^n$ is (for some/all $n\in\Z\setminus \{0\}$) and if and
only if some conjugate of $g$ is.

We now examine which types of isometry are preserved when moving from
$G$ acting on an arbitrary metric space $(X,d_X)$ to the action of $G$ on a
Rips orbit graph $\Gamma_r(x_0)$ (assumed connected with path metric
$d_\Gamma$). The function $f$ going the other way
from $\Gamma_r(x_0)$ to $X$ (say with edges
of $\Gamma_r(x_0)$ first sent to a nearest vertex) is a coarse
$G$-equivariant map satisfying $d_X(f(v),f(w))\leq rd_\Gamma(v,w)$ for
vertices $v,w$ and hence $d_X(f(x),f(y))\leq (r+1) d_\Gamma(x,y)$ for
$x,y$ arbitrary points.
If this were also a quasi-isometric 
embedding then the three types of isometry would be the same in both actions
of $G$. However we have already seen that this will often not be the case.
In complete generality,
we have that $g\in G$ being a loxodromic element is preserved because
$d_X(v,w)\leq r d_\Gamma(v,w)$ for vertices $v,w$.
By Proposition \ref{pres} we know that
$G$ having unbounded orbits is preserved, so applying this to 
$\langle g\rangle$ tells us that parabolics go to parabolics or loxodromics.
But if we take the parabolic element $f(z)=z+1$ acting on the upper
half plane $\HH^2$ with the hyperbolic metric and $x_0=i$ then $\Gamma_1(i)$
is the simplicial line and now $f$ acts loxodromically. Note that this action
on $\HH^2$ has locally finite orbits and both spaces are hyperbolic.

However in general elliptic elements can go to any of elliptic, parabolic
or loxodromic: for instance our example in Proposition \ref{pres} had
$X$ bounded, so all elements of $G$ act elliptically, but the Rips orbit
graph returned the Cayley graph of $G$ where all three types of isometry
are possible.

If we now assume though that the orbit $G\cdot x_0$ is locally finite
then Proposition \ref{pres} also told us that having bounded orbits
is preserved. Thus here elliptics and loxodromics will remain as such
in their action on the Rips orbit graph, so the only ambiguity is that
elements acting parabolically on $X$ might act parabolically or
loxodromically on $\Gamma_r(x_0)$. In fact we mention here that if both $X$ and
$\Gamma_r(x_0)$ are hyperbolic spaces with $G\cdot x_0$ a locally finite
orbit then parabolic isometries in $X$ will always be loxodromic
in $\Gamma_r(x_0)$ (this will be seen in Section 8). 

Let us now concentrate on the case where $X$ is a geodesic metric space
satisfying any of the equivalent definitions of $\delta$ - hyperbolicity.
We will refer to this throughout as $X$ is hyperbolic but no further
conditions such as properness of $X$ will be assumed. 
As $X$ is hyperbolic we can look at the
action of $G$ by homeomorphisms (though not isometries)
on the (Gromov) boundary $\partial X$ of $X$ to obtain the limit
set $\partial_G X$ which is a subset of $\partial X$. This subset is
$G$-invariant and we have $\partial_H X\subseteq \partial_G X$ if $H$ is
a subgroup of $G$. In particular we can take an arbitrary isometry $g$
of some hyperbolic space $X$ and consider $\partial_{\langle g\rangle} X$.
We can summarise the facts we will need in the following well known proposition.
\begin{prop} \label{cyc}
(i) If $g$ is elliptic then $\partial_{\langle g\rangle} X=\emptyset$. This
is always the case if $g$ has finite order but might also occur for elements
with infinite order. In either case $g$ might fix many or no points on
$\partial X$.\\
(ii) If $g$ is loxodromic then $\partial_{\langle g\rangle} X$ consists
of exactly 2 points $\{g^\pm\}$ for any $x\in X$ 
and this is the fixed point set of $g$ (and $g^n$
for $n\neq 0$) on $\partial X$. Moreover for any $x\in X$ we have
$g^n(x)\rightarrow g^+$ and $g^{-n}(x)\rightarrow g^-$ in $X\cup\partial X$
as $n\rightarrow\infty$.\\
(iii) If $g$ is parabolic then $\partial_{\langle g\rangle} X$ consists
of exactly 1 point and again this is the fixed point set of $g$ (and $g^n$
for $n\neq 0$) on $\partial X$.\\
\end{prop}

Moving back now from cyclic to arbitrary groups, for any group $G$ 
(not necessarily finitely generated)
acting by isometries on an arbitrary hyperbolic space $X$, we have the
Gromov classification dividing possible actions into five very different
classes. In the case that $G=\langle g\rangle$ is cyclic, the first three 
classes correspond to the three cases in Proposition \ref{cyc} and
the next two do not occur (for these facts and related references, see 
\cite{abos} and \cite{ren}):\\
\hfill\\
(1) The action has {\bf bounded orbits}. This happens exactly when
$\partial_G X$ is empty, in which case all elements are elliptic.
However if all elements are elliptic then we can also be in case (2).\\
\hfill\\
(2) The action is {\bf parabolic} (or horocyclic),
meaning that $\partial_G X$ has exactly
one point $p$. Note that, despite the name, we
can still be in this case but with $G$ consisting only of elliptic elements
(there exist both finitely generated and non finitely generated examples).
However any non elliptic element must be parabolic with limit
set $\{p\}$. Such an action cannot be cobounded.\\
\hfill\\
(3) The action is {\bf lineal}, meaning that
$\partial_G X=\{p,q\}$ has exactly 2 points (we use dihedral to indicate
that these points can be swapped and linelike to say that they are pointwise
fixed, thus if the action of $G$ is dihedral then there is an index 2
subgroup $H$ of $G$ with the action of $H$ being linelike). In this case there
will exist some loxodromic element in $G$ with limit set $\{p,q\}$ and
indeed all loxodromics have this limit set. Moreover $G$ contains
no parabolics (such an element $g$ would have to fix $p$ or $q$ but if
it were $p$ then $g$ would move $q$ outside $\{p,q\}$ which is
$G$-invariant) but there might be many elliptic elements in $G$.\\
\hfill\\
(4) The action is {\bf quasi-parabolic} (or focal).
This says that the limit set
has at least 3 points, so is infinite, but there is some point
$p\in\partial_G X$ which is globally fixed by $G$. This implies that
$G$ contains a pair of loxodromic elements with limit sets
$\{p,q\}$ and $\{p,r\}$
for $p,q,r$ distinct points.
Again despite the name, the group $G$ need not
contain parabolic elements but any that do exist will have limit set
$\{p\}$ (otherwise they would move $p$) and again there can be many elliptic
elements in $G$.\\
\hfill\\
(5) The action is {\bf general}: the limit set is infinite and we have two
loxodromic elements with disjoint limit sets, thus by ping pong high powers
of these elements will generate a quasi-embedded copy of the free group
$F_2$. In this case $G\setminus\{e\}$
might only consist of loxodromic elements,
or contain just parabolics as well, or contain just elliptics as well, or
all three types of isometry.

We refer to these descriptions as the hyperbolic type of the action.
We also call the actions from (1) to (4) {\bf elementary},
so that any isometric action on a hyperbolic space is either of general
type or is elementary. Whilst there is agreement that types (1) to (3)
are deserving of the title ``elementary'', there are different conventions
on whether or not to include type (4). Here we will regard type (4)
as elementary because we will be interested in versions of the Tits
alternative and type (5) is the only one that guarantees a non abelian
free subgroup.

We now give some examples that illustrate this classification and which
will be relevant later. For type (1), any finite group must act in this way
and any group can act trivially on any hyperbolic space. A rather awkward
variation on the latter is to take any infinite group $G$ and the complete
graph $\Gamma$ with vertex set the elements of $G$. Then this graph is
bounded and hence a hyperbolic space, with $G$ acting freely (if order two
elements are taken care of) and even topologically properly.

Moreover any
infinite group acts by isometries on its Groves - Manning combinatorial
horoball (see \cite{mn3} 3.1) which is a hyperbolic space. This is a parabolic
action with the finite order elements acting elliptically (as always)
and all infinite order elements acting parabolically. Moreover if $G$
is any infinite group which is finitely generated then this
construction results in a locally finite (but never bounded valence: see
Section 8)
hyperbolic graph with $G$ acting metrically properly and freely.
This also means that any countable group has such an action by embedding
it in a finitely generated group.

Actions of type (5) can be obtained by taking a non elementary
word hyperbolic group acting (metrically properly and coboundedly) on its
Cayley graph $C$ with respect to a finite generating set, so that $C$
has bounded valence. If we take an arbitrary subgroup $H$ of $G$ and restrict
this action to $H$ then it is no longer cobounded, unless $H$ has finite
index in $G$, but it will be of type (5), or of type (3) if $H$ is infinite
and virtually cyclic, or of type (1) if $H$ is finite.

We might wonder which of these actions occur if we restrict our spaces to
simplicial trees, whereupon we can certainly get types (1), (3) and (5).
It might be
expected that no action of type (2) exists on a tree because there are
no parabolic isometries. However, given any countably infinite group $G$
which is locally finite (every finitely generated subgroup is finite),
we can obtain an action of $G$ on a locally finite tree with unbounded orbits
where every element acts elliptically. This is sometimes called the
coset construction and is described in \cite{ser} I.6.1 or in  \cite{bh}
II.7.11.
The basic idea is to express $G$ as an increasing union
$\cup_{i=0}^\infty G_i$ of finite subgroups $G_i$ with $G_0=\{e\}$. The 0th
level vertices are the elements of $G$, the 1st level vertices are the
cosets of $G_1$ and so on. The stabiliser of a vertex in level $i$
has size $|G_i|$ and its valence is
$[G_i:G_{i-1}]+1$ (and 1 in level zero). Thus this action is metrically
proper and we even have a bounded valence tree
if $[G_i:G_{i-1}]$ is bounded. However we note that here $G$ is certainly
not finitely generated, and indeed if any finitely generated group acts
on a tree with all elements elliptic then it has a global fixed point by
\cite{ser} I.6.5: moreover for a generating set $\{g_1,\ldots ,g_n\}$ it
is enough to check that each $g_i$ and $g_ig_j$ is elliptic.

This just leaves actions of type (4), which can occur for finitely
generated groups on bounded valence trees. One example is
the action of the Baumslag - Solitar
group $BS(1,2)$ on its Bass - Serre tree using the standard HNN extension
decomposition. Another is to take the wreath product $C_p \wr \Z$ for
$p$ a prime and consider the action on its Bruhat - Tits tree.

In these two examples, it might be noted that we have a homomorphism from
our group $G$ to $\Z$ where the elements outside the kernel are exactly the
elements acting loxodromically.
We might hope that the same result holds for all actions of type (3) and
(4), or at least (thinking of $\R$ acting on itself by translations) a
homomorphism from $G$ to $\R$ with the same property.
But in general we have to use homogeneous
quasi-morphisms $q:G\rightarrow \R$ rather than homomorphisms.
We have by \cite{mn3} Proposition
4.9 that if $G$ acts on $X$ and fixes a point in $\partial X$ (so this
might be an action
of type (1) or (2), but then the conclusion is empty, or type (3) on possibly
dropping to an index 2 subgroup, or of type (4)) then by use of Busemann
functions there exists a homogeneous
quasi-morphism $q$ (sometimes called the Busemann quasicharacter)
where the quasi-kernel ($\{g\in G\,|\,q(g)=0\}$) consists
exactly of the elements of $G$ which do not act loxodromically.
Sometimes this will be a genuine homomorphism, for instance if $X$ is
a tree and moreover if the action is by tree automorphisms (so ruling
out $\R$ and its subgroups acting on itself by translation) then this
homomorphism can be taken to have image in $\Z$. Other cases are in
\cite{cdcmt} where it is shown that the Busemann quasicharacter is also
a genuine homomorphism to $\R$ if $G$ is amenable or $X$ is proper.
In particular we might
expect to find examples of groups with very restricted actions on hyperbolic
spaces amongst those with no unbounded quasi-morphisms to $\R$. Indeed it
is shown in \cite{mn3} that $SL(n,\Z)$ for $n\geq 3$ (amongst other groups)
only has actions of types (1) and (2); in particular there is no action of
$SL(n,\Z)$ on any hyperbolic space where any element is loxodromic
(see also \cite{clvg} for the specific case of trees and \cite{hae} for
a more general result).

Suppose that $G$ has two isometric actions on the hyperbolic spaces
$Y$ and $X$ and that we have a coarse $G$-equivariant quasi-isometric
embedding $F$ from $Y$ to $X$. We have already said that this preserves
the isometry type of individual elements and indeed $F$ also preserves
the hyperbolic type of the action. This is shown in \cite{abos} Lemma 4.4
when $F$ is also a quasi-isometry but the proof goes through for this case too.

As our main aim in this
paper is to obtain cobounded actions from arbitrary actions using the Rips
orbit graph, we finish this section with a definition which, as
noted in \cite{abos},
is a good way to obtain cobounded actions from non-cobounded ones
when our space $X$ is hyperbolic.

\begin{defn} A subset $S$ of a hyperbolic space $X$ is
$K$-quasiconvex if for any $x,y\in S$ and any geodesic segment
$[x,y]$ in $X$, every point of $[x,y]$ is within distance $K$
of a point of $S$.
\end{defn}

\begin{prop} \label{os}
(\cite{abos} Proposition 3.14)\\
Let $G$ act on a hyperbolic space $X$ by isometries and with a $K$-quasiconvex
orbit $G\cdot x_0$ for some $x_0\in X$. Then the Rips orbit graph 
$\Gamma_{2K+1}(x_0)$ is connected and the path metric $d_\Gamma$ satisfies
\[d_\Gamma(a,b)\leq d_X(a,b)\leq (2K+1)d_\Gamma(a,b).\]
Thus the map $\phi:\Gamma_{2K+1}(x_0)\rightarrow X$ which is the identity on 
$G\cdot x_0$ and sends an edge of $\Gamma$ to a geodesic in $X$ between
its endpoints is a coarse $G$-equivariant quasi-isometric embedding.
\end{prop} 

It is the case that in (4), where we have a quasi-parabolic action of a
group on a hyperbolic metric space, all orbits are quasi-convex.

\section{Group actions with two limit points}

We first mention some standard definitions and results.

Given a metric space $X$ and constants $K\geqslant 1,C\geqslant 0$,
a $(K,C)$-(discrete) {\bf quasigeodesic} is a $(K,C)$ quasi-isometric
embedding from $\mathbb R$ ($\Z$) to $X$. 
On being given any isometry $g$ and any point $x_0$ of an arbitrary metric
space $X$, we have that the orbit of $x_0$ under $\langle g\rangle$ is a
discrete quasigeodesic if and only if $g$ is a loxodromic element.
Indeed if $g$ is loxodromic then we can take $C=0$ since we saw in Section 3
that
$\tau(g)|m-n| \leqslant d(g^m(x_0),g^n(x_0))$. Thus the lower bound is
independent of $x_0$, but we also have
\[d(g^m(x_0),g^n(x_0))=d(g^{m-n}(x_0),x_0))\leqslant |m-n| d(g(x_0),x_0)\]
so that the upper bound will generally depend on $x_0$.

Suppose now that $X$ is any geodesic metric space (here geodesics need not be
unique) and we have an isometry $g$ of $X$ which is acting loxodromically.
In general we cannot assume that there is a (two way infinite) geodesic
of $X$ which is invariant under the action of $\langle g\rangle$. However
we can obtain a continuous quasigeodesic from the orbit of $x_0$ by first
joining $x_0$ to $g(x_0)$ by some geodesic segment $\gamma$ of length
$l:=d(g(x_0),x_0)>0$ and then for each $n$ 
in $\Z\setminus\{0\}$ we join $g^n(x_0)$ 
to $g^{n+1}(x_0)$ by the image of $\gamma$ under $g^n$, which is also
a geodesic segment of length $l$.
We can then define $f:\R\rightarrow X$ by sending $[n,n+1]$ evenly to
this geodesic $g^n(\gamma)$, so we first map $[n,n+1]$ linearly to
$[ln,l(n+1)]$ and then map this isometrically to $g^n(\gamma)$.
Now suppose that we have $s,t\in\R$ with $s\in [m,m+1]$ and $t\in [n,n+1]$.
We can then see that
\[\tau(g)|m-n|-2l\leq d(f(s),f(t))\leq l|m-n|+2l\]
but as we have $|t-s|-1\leq |m-n|\leq |t-s|+1$, we obtain
\[\tau(g)|s-t|-\tau(g)-2l\leq d(f(s),f(t))\leq l|s-t|+3l.\]
In particular on setting $k=\mbox{max}(1/\tau(g),l)$
we have a $(k,3k)$-quasigeodesic
containing the orbit of $x_0$ under $\langle g\rangle$ and whose image is
$\langle g\rangle$-invariant.

So far we have not used anything about $X$ other than it is a geodesic
metric space. But let us now say that $X$ is a $\delta$-hyperbolic space
and $g$ is a loxodromic isometry of $X$ with fixed points $g^\pm$ on
$\partial X$. We know that on
$X\cup\partial X$, it is the case that for any $x\in X$ we have
$g^n(x)\rightarrow g^+$ as $n\rightarrow\infty$ and $g^n(x)\rightarrow g^-$ as
$n\rightarrow -\infty$, so our $(k,3k)$-quasigeodesic has
endpoints $g^\pm$. We now utilise the stability of quasigeodesics
with the same end points, by quoting a version of this property
that is Lemma 2.11 in \cite{man}.
\begin{lem} \label{mann}
Let $K\geq 1$, $C\geq 0$ and $\delta\geq 0$. Then there is some
$B=B(K,C,\delta)$ such that if $X$ is any $\delta$-hyperbolic geodesic metric 
space and $\gamma$ and $\gamma'$ are two
$(K,C)$-quasigeodesics in $X$ with the same endpoints in $X\cup\partial X$
then the image of $\gamma$ lies in a $B$-neighbourhood of the image of 
$\gamma'$.
\end{lem}

Let us now consider the case where an arbitrary group $G$ has an action
of type (3) on an arbitrary hyperbolic space; that is
$\partial_GX$ has two limit points.
Thus there will be some $g\in G$ acting loxodromically with these two limit
points equal to its fixed points on the (Gromov) boundary $\partial X$.
We have mentioned that a necessary condition for the existence of such an
action, at least on possibly dropping down to an index 2 subgroup if this
action is dihedral,
is a homogeneous quasi-homomorphism $q:G\rightarrow \R$ with
$q(g)\neq 0$. Then the quasi-kernel consists of those elements acting
elliptically and everything else acts loxodromically. As pointed out in
\cite{abos} subsection 4.3, we have a converse: given such a $q$ which is
non zero, we have $C>0$ such that the (usually infinite)
subset $S=\{g\in G\,|\,|q(g)|<C\}$ generates $G$ and the Cayley graph
of $G$ with respect to $S$ is quasi-isometric to $\R$, thus $G$ has an
action of type (3) on this Cayley graph. In particular, on taking a
torsion free
hyperbolic group $H$ which does not have a surjective homomorphism to $\Z$
we have by \cite{epfj}
quasi-morphisms of $H$ which give rise to actions of type (3)
that are linelike but which are not witnessed
by homomorphisms from $H$ to $\R$.   

Now suppose we have a finitely generated group $G$ with a type (3)
action and a genuine homomorphism $r$ from $G$ to $\R$ where the kernel
consists of all elements acting elliptically. Although for any $g\in G$
which is loxodromic under this action, we can find using $r$ a homomorphism
$\theta$ to $\Z$ with $\theta(g)\neq 0$, we might not be able to do this
simultaneously for all loxodromic elements. For instance $\Z^2$ acting
on $\R$ with generators equal to the maps $x\mapsto x+1,x\mapsto\sqrt{2}$
has all non identity elements acting loxodromically. We can even have this
situation when acting by automorphisms on a hyperbolic graph, for instance
on applying the original Rips construction to this action  on $\R$, we get
an action of $\Z^2$ by automorphisms on a hyperbolic graph where every
non identity element acts loxodromically and where the action is free
and topologically proper.

Having considered this range of examples, we notice that none of these
actions have a locally finite orbit. With just this condition in place, we
are able
to get a full result on identifying type (3) actions via homomorphisms to
$\Z$. Indeed for our next theorem, no
assumptions on $G$, $X$ or the action of $G$ on $X$ will be assumed other
than those which are expressly stated. In particular there is
no assumption here of $G$ of being finitely generated or even countable:
even in the case where
$X$ is a locally finite graph and $G$ acts by automorphisms, it can happen
that $G$ is uncountable but still acts faithfully (see below).

\begin{thm} \label{two}
Suppose that the group $G$ acts by isometries on the geodesic 
$\delta$-hyperbolic space $X$ with a locally finite orbit.
If the limit set $\partial_GX$ consists of exactly two points
(so that all elements of $G$ are either elliptic or loxodromic)
then $G$ is of the form $E\rtimes\Z$ or $E\rtimes(C_2*C_2)$,
according to whether each limit point is fixed by $G$ or not,
where $E$ consists entirely of elliptic elements
and acts on $X$ with bounded orbits. In the first case the set of elliptic
elements of $G$ is exactly $E$. In the second
case the set of elliptic elements is exactly $E$ union the complement of
the natural index 2 subgroup $E\rtimes\Z$.
\end{thm}
\begin{proof}
As $G$ fixes the
limit set $\{x^+,x^-\}$ say, we have that $G$ or an index two subgroup 
(renamed $G$ from now until the end of the proof) fixes both $x^+$ and 
$x^-$ pointwise. 
In $G$ there are loxodromic elements with attracting fixed point $x^+$
and repelling fixed point $x^-$. Take any one of those and call it $h$.

We will now build a related graph $\Gamma_0$ to work with by following
the construction of the Rips orbit graph in Definition \ref{rog}, though
we cannot assume that $G$ is finitely generated.
First take 
this loxodromic element $h\in G$ and use
it and the point $x_0\in X$ which has a locally finite orbit under $G$
to form a $(k,3k)$ quasi-geodesic 
$\delta_0$ as above.
We let our set $V(\Gamma_0)$ of vertices to be the orbit of $x_0$ under $G$
which is of course $G$-invariant. Now any point $g(x_0)$ in this orbit
will lie on the $(k,3k)$ quasi-geodesic $g(\delta_0)$ which has the same
endpoints as $\delta_0$. Thus on setting $B$ to be the Morse
constant $B(k,3k,\delta)$ we have that for all $g(x_0)\in G\cdot x_0$ 
there exists a point $p$ on $\delta_0$ with distance at most $B$ from 
$g(x_0)$, thus there exists $m\in\Z$ such that 
$d(g(x_0),h^m(x_0))\leq K:=B+l/2$. As this orbit is locally
finite, there is a finite list $g_1(x_0)=x_0,\ldots ,g_s(x_0)$ of all the 
points in $G\cdot x_0$ which are at most distance $K$ from $x_0$. Hence 
as $d(h^{-m}g(x_0),x_0)\leq K$ we have 
$h^{-m}g(x_0)=g_i(x_0)$ for some $i$ and so $g(x_0)=h^m(g_i(x_0))$, showing us
that $h$ has finitely many orbits in its action on $G\cdot x_0$. 

We now shrink the list $g_1(x_0),\ldots ,g_s(x_0)$ so that it contains
exactly one point from each $\langle h\rangle$ orbit. This now means that
every point in $G\cdot x_0$ has a unique representation in the form
$h^mg_i(x_0)$, because if $h^mg_i(x_0)=h^ng_j(x_0)$ then $i=j$ but 
$g_i^{-1}h^{n-m}g_i$ is loxodromic if $n\neq m$ and so has no fixed points. 
In order to turn
our collection of vertices $G\cdot x_0$ into our graph $\Gamma_0$, we
set $M$ to be the maximum of $d(hg_i(x_0),g_i(x_0))$ for $i=1$ to $s$.
We now take the Rips orbit graph $\Gamma_r(x_0)$ where $r$ is 
the maximum of $K$ and $M$.
We then define $\Gamma_0$ to be $\Gamma_r(x_0)$ for this value of $r$,
with $G$ acting on $\Gamma_0$ by graph automorphisms and
transitively on its vertices. Also
$\langle h\rangle$ has finitely many orbits of vertices
and $\Gamma_0$ is a locally finite graph (as the orbit is locally
finite) and so of bounded valence.
Moreover $\Gamma_0$ is connected. This is because any vertex in $\Gamma_0$
is of the form $h^mg_i(x_0)$ and thus is at distance at most $M$ in
$X$ from $h^{m+1}g_i(x_0)$ and $h^{m-1}g_i(x_0)$. Thus we can go
up and down this $\langle h\rangle$ orbit by a path in $\Gamma_0$ until we
reach $g_i(x_0)$, but then this is at most $K$ in $X$ from
$x_0$, thus we can reach $x_0$ in $\Gamma_0$ as well. 

We also have quasiconvexity of the orbit $G\cdot x_0$: indeed here we
can show directly that there is $\lambda>0$ and $\epsilon$ with
\[\lambda d_{\Gamma_0}(v,w)-\epsilon
  \leq d_X(v,w)\qquad\mbox{for }v,w\in G\cdot x_0.\]
As we can take $v=h^mg_i(x_0)$ and $w=h^ng_j(x_0)$ with $k=n-m>0$
without loss of generality, we obtain
\begin{eqnarray*}
  d_X(h^mg_i(x_0),h^ng_j(x_0))&=&d_X(g_i(x_0),h^kg_j(x_0))\\
                              \geq
d_X(h^kg_j(x_0),g_j(x_0))-d_X(g_j(x_0),g_i(x_0))&\geq& k\tau(h)-2K
\end{eqnarray*}
whereas from the above we know that in $\Gamma_0$ we have
$d_{\Gamma_0}(h^mg_i(x_0),h^ng_j(x_0))=d(g_i(x_0),h^kg_j(x_0))$
but we can get from $g_i(x_0)$ to $g_j(x_0)$ in at most 2 steps and
then from $g_j(x_0)$ to $h^kg_j(x_0)$ in at most $k$
so we can take $\lambda=\tau(h)>0$ and $\epsilon=2(\tau(h)+K)$.

Thus we now have that $G$ acts on the graph $\Gamma_0$ vertex
transitively by automorphisms and indeed coboundedly by isometries.
Moreover $\Gamma_0$ is a hyperbolic graph as it quasi-isometrically embeds
in $X$ and the action of $G$ on $\Gamma_0$ is also lineal, so that
$\partial_G\Gamma_0$ consists of two points which we also
call $x^+,x^-$. 
We can now finish off our result by purely combinatorial means,
on invoking results in \cite{srv} and others about vertex transitive,
connected, locally finite graphs with polynomial (here linear) growth.
These results do not mention hyperbolicity, although underlying them
is Gromov's famous theorem on groups of polynomial growth.

Alternatively we can proceed by invoking results on locally compact
topological groups as in \cite{cdcmt}. If $\Gamma_0$ is equal to the
simplicial line then we are done. Otherwise
we replace $G$ with the
automorphism group $A$ of $\Gamma_0$, which is equal to its isometry group.
Note that if we have established
that $A=E\rtimes\Z$ as in the statement of the theorem,
corresponding to when there is no isometry swapping $x^+$ and $x^-$,
then any subgroup
of $A$ containing a loxodromic element will not lie in $E$, so
we can take $G$, or at least its image in $A$ if this is not faithful
(whereupon $G$ is still of the same form),
equal to $(E\cap G)\rtimes\langle t\rangle$. Here $E\cap G$ also
acts on $\Gamma_0$ (and hence on $X$) with bounded orbits and $t$ is
any element of $G$ mapping to an integer generating the image
of $G$ in $\Z$. 

If however we have $A=E\rtimes (C_2*C_2)$ then we have elements swapping
$x^+$ and $x^-$, with the index two subgroup $B=E\rtimes\Z$
fixing those two points pointwise. So any element $g$
in $A\setminus B$ swaps $x^+$ and $x^-$ on the boundary $\partial X$.
and therefore must be elliptic:
if $g$ were loxodromic then $g^2$ is also loxodromic but it would have 
at least four fixed points on the boundary. Thus under the natural
homomorphism from $A$ to $C_2*C_2$, any subgroup $G$ containing a loxodromic
element will have infinite image, which will thus be equal to $\Z$ or
$C_2*C_2$. The first case is as before and the second case follows by
the same argument.

Thus we put the compact open topology (or equivalently the pointwise
topology) on $A=Isom(\Gamma_0)=Aut(\Gamma_0)$. But $\Gamma_0$ is of course
a proper metric space so $A$ is a second - countable and locally compact group,
and the action of $A$ on $\Gamma_0$ is continuous and proper in the topological
group sense by \cite{dcdlh} Lemma 5.B.4. Moreover as the action of $A$ on
$\Gamma_0$ is cocompact, we have that $A$ is compactly generated.
Now \cite{cdcmt} Proposition 5.6 applies to locally compact groups with a
continuous proper cocompact isometric action on a proper hyperbolic space
$X$ with two limit points, thus to $A$ acting on $\Gamma_0$. It states
that $A$ has a maximal compact normal subgroup $E$ so that $A/E$ is
isomorphic to $\Z$, $\R$, $C_2*C_2$ or $Isom(\R)$ where the second two
cases have an index 2 subgroup corresponding to the first two cases
respectively
and the homomorphism to $\Z$ or $\R$ is the Busemann quasicharacter
at one of the limit points.

Now a compact subgroup of $A$ would have bounded orbits in $\Gamma_0$
because the action of $A$ is proper. Moreover in our case we must
have that $A/E$ is countable: if not then take representatives
$\{a_i\in A\,|\,i\in I\}$ for the cosets of $E$ in $A$ and define a map
from this uncountable
set to the vertices of $\Gamma_0$ via $a_i\mapsto a_i(x_0)$.
Then we will have $a_i(x_0)=a_j(x_0)$ for some $i\neq j$ and so the
element $a_i^{-1}a_j$ fixes $x_0$. Thus it is elliptic and so is sent
to zero by the Busemann homomorphism, so lies in $E$ which is a contradiction.
Thus $A$, and hence $G$, is of the form $E\rtimes\Z$ or $E\rtimes(C_2*C_2)$
as in the statement of the theorem.
\end{proof}

{\bf Example}: To show that $G$ can be uncountable, consider 
the graph consisting of two copies of the standard Cayley graph of $\Z$,
with vertices labelled $(n,1)$ and $(n,2)$ in each copy, and with further
edges where $(n,1)$ is joined to $(n-1,2)$ and $(n+1,2)$. Then for a given
$n\in\Z$ we have that $(n,1)$ and $(n,2)$ have the same neighbours so
permuting these two vertices and fixing every other vertex is a graph
automorphism, and this can be done independently of any other $n$. Thus
the resulting group of automorphisms so obtained is 
$\left(\Pi_{n\in\Z} C_2\right) \rtimes \Z$ where $\Pi_{n\in\Z} C_2$ 
is the Cartesian (unrestricted) direct product.
 
We can
put Theorem \ref{two} in the context of coarse $G$-equivariant
quasi-isometric embeddings which we are using to compare actions.
Note that a surjective homomorphism from $G$ to $\Z$ or $C_2*C_2$ is a
necessary condition for such a map below to exist.

\begin{co} \label{cgqie}
Suppose that an arbitrary group $G$ acts by isometries
  on an arbitrary geodesic $\delta$-hyperbolic space $X$. If the
  limit set $\partial_GX$ consists of two points and there is a
  locally finite orbit, as in Theorem \ref{two}, then we have an
  action of $G$ by automorphisms on the simplicial line $L$ and a coarse
  $G$-equivariant quasi-isometric embedding from $L$ to $X$.
\end{co}

\begin{proof}
  If we are in the first case of Theorem \ref{two} then let $h\in G$ be
  a loxodromic element generating the $\Z$ quotient. If we are in the
  second case, then drop to the relevant index two subgroup to get this
  element $h$ and let $t$ be any element in $G$ generating the $C_2*C_2$
  quotient along with $h$: note that $th^it^{-1}=h^{-i}$ in this quotient.

  Our action of $G$ on the simplicial line $L$ is defined by making the
  the quotient $\Z$, which is generated by $h$,  act as $h^i(l)=l+i$ and
  $t$ act as $t(l)=-l$, with the kernel $E$ acting trivially.
  Our quasi-isometric embedding $F:L\mapsto X$ (or from $L$ to $\Gamma_0$
  and then composed with the quasi-isometric embedding established in the
  proof from $\Gamma_0$ to $X$; the details are the same)
  is given by $F(m)=h^m(x_0)$
  for $m\in\Z$ and the interval $[m,m+1]$ is mapped to a geodesic
  between $h^m(x_0)$ and $h^{m+1}(x_0)$ as in the start of the proof
  of Theorem \ref{two} (the original element $h$ might not generate
  the $\Z$ quotient but some loxodromic element will a fortiori and the proof
  of Theorem \ref{two} applies equally for this element too).
We saw that $F$ is
a quasi-geodesic, thus a quasi-isometric embedding from $L$ to $X$
(or $\Gamma_0$).
  
  To show the coarse $G$-equivariance, let us first do this for an integer
  $m\in L$. As Theorem \ref{two} tells us that $E$ acts with bounded
  orbits (on $X$ or on $\Gamma_0$),
  we have $L>0$ such that for all $k\in E$ both
  $d(x_0,k(x_0))\leq L$ and $d(x_0,tk(x_0))\leq L$ hold, where
the second inequality follows because $d(x_0,tk(x_0))\leq
  d(x_0,t(x_0))+d(t(x_0),tk(x_0))$.
  
  First take any element of the form $h^ik\in G$, where $i\in\Z$ and
  $k\in E$. Then $d(F(h^ik(m)),h^ik(F(m)))=d(h^{i+m}(x_0),h^ikh^m(x_0)$
  but $kh^m=h^mk'$ for $k'$ also in $E$, so this is equal to
  $d(h^{i+m}(x_0),h^{i+m}(k'(x_0)))=d(x_0,k'(x_0))\leq L$. That covers all
  group elements in the first case, whereas in the second we also need
  to check elements of the form $th^ik$. The same calculation gets us
  to $d(h^{-i-m}(x_0),th^ikh^m(x_0))$ but $th^ikh^m(x_0)=h^{-i-m}tk'(x_0)$
  thus this distance is at most $L$ too.

  Finally for a point $l$ on $L$ lying between $m$ and $m+1$, we have
  for any $g\in G$ that $d(F(g(l)),g(F(l)))$ is at most
 \[ d(F(g(l)),F(g(m)))+d(F(g(m)),g(F(m)))+
  d(g(F(m)),g(F(l))),\]
  with the second term at most $L$ from above
  and the first and third terms bounded
  independently of $l$ and $g$ because $F$ is a quasi-isometric embedding
  and $g$ acts as an isometry on both $L$ and $X$.
  \end{proof}

\section{Applications to quasitrees}

A quasitree is a geodesic metric space
which is quasi-isometric to a tree. (Here a quasitree need not be a
combinatorial graph, though many of the examples that interest us will be of
that form. However a tree will always refer to a simplicial tree, that is
a connected combinatorial graph, equipped with its path metric, that contains
no closed loops.) In the literature are various results about how groups act
on trees. Perhaps the most famous (from \cite{ser}) is that a group acting
freely by automorphisms on a tree is a free group: there is no finiteness
assumption here (either on cardinality or generation of the group, nor
on local finiteness of the tree). If we now change acting freely by
automorphisms to acting topologically properly
by automorphisms, that is every
vertex has a finite stabiliser, then for a group $G$ acting in this way on
an arbitrary tree it is the case that $G$ is virtually
free if $G$ is finitely generated, but not in general which we have seen
with coset constructions (even for bounded valence trees).

We also have (again in \cite{ser}) that if a finitely generated group
$G$ acts on a tree without a global fixed point, which is equivalent to
acting with unbounded orbits, then $G$ splits as an HNN extension or
amalgamated free product.

We might now ask: suppose a group $G$ acts on some quasitree in one of the
ways given above then can we draw the same conclusion about the structure
of $G$? A moment's thought reveals many counterexamples:\\
$\bullet$\,\,
If $G$ is any group whatsoever acting on its complete graph $\Gamma$ by
automorphisms then this action is free and topologically proper, with
$\Gamma$ a quasitree because it is quasi-isometric to a point. (However
$\Gamma$ is not locally finite and the action has bounded orbits and is
not metrically proper.)\\
$\bullet$\,\,
Even for geometric and free actions of a finitely generated group $G$ on a
bounded valence graph which is a quasitree, $G$ need not be free but could be
virtually free,
eg the Cayley graph of a virtually free group will be a quasitree. \\
$\bullet$\,\,
Even if the quasitree is a bounded valence graph then we have mentioned the
case of a subgroup of $\R$ acting on $\R$ freely by isometries but
not by automorphisms. We have also pointed out after Lemma \ref{mann} that
this can be turned into a free and topologically proper (but not metrically
proper) action by automorphisms on an infinite valence graph quasi-isometric
to $\R$ by using the Rips construction.
Alternatively we can think of these actions on $\R$ as being free actions
by isometries on an $\R$-tree (and any $\R$-tree is a quasitree; for instance
see \cite{alc}). Here we also have free actions of other groups on $\R$-trees
too, such as surface groups (which can then be turned similarly into free
actions by automorphisms on graphs which are quasitrees).

Thus free actions of groups on quasitrees appear to have little to do with
free actions of groups on trees. If we now move to proper actions, the coset
construction gives us proper actions of infinitely generated groups on
bounded valence trees but these groups are locally finite and so are a long
way from virtually free groups. Perhaps the best we can hope for is to say
something about a finitely generated group acting metrically properly
on an arbitrary quasitree.

Turning now away from proper and/or free actions to unbounded actions of
groups on quasitrees, the situation seems even more extreme. Given a
group $G$ with an unbounded quasi-morphism, we mentioned after Lemma
\ref{mann} the construction in \cite{abos} which gives us an unbounded action
of $G$ by automorphisms on a (generally infinite valence) graph which
is a two ended quasitree. Thus a hyperbolic group with property (T) has
an unbounded orbit on a quasitree but every action of any finite index
subgroup on any tree has a global fixed point.

There is an even greater contrast in \cite{balas} (as also
mentioned in Section 7). It is shown there that any acylindrically hyperbolic
group has an action on a quasitree which is (acylindrical and)
not elementary. Thus it seems that vast classes of groups have some
unbounded action on a quasitree even if all their actions on trees have
a global fixed point. Indeed even results giving us sufficient conditions
for a given action of a group on a tree to have bounded orbits do
not have an equivalent for quasitrees. For instance a well known consequence
of \cite{ser} I.6.5 is that if a group $G=\langle x,y\rangle$ acts on a tree
with $x$, $y$ and $xy$ all elliptic elements then the action has bounded
orbits. This also holds for $\R$-trees by \cite{mrsh}.
But in \cite{mn1} the example is given of the Farey tessellation
of the hyperbolic plane, which becomes a combinatorial graph when given
the path metric. This graph satisfies the bottleneck criterion mentioned
below and so is a quasitree.
Now the modular group $SL(2,\Z)$ acts on the hyperbolic
plane with $x$ an element of order 2, $y$ an element of order 3 and their
product a parabolic element of infinite order. This fixes a rational
point on the boundary and hence acts elliptically on this graph, but
the overall action is unbounded. Even more extreme examples come from
the result in \cite{balas} mentioned above, because we can take a hyperbolic
group generated by $x_1,\ldots ,x_n$ where all group words up to a given length
have finite order, thus must act elliptically, but we still have an unbounded
action on a quasitree. (However note that these are graphs with infinite
valence and also that we have not given an example
of an unbounded action of a finitely generated group on any quasitree
where every element is elliptic.)
Finally in our list of examples, we do have groups which can
only act on a quasitree with bounded orbits, such as $SL(n,\Z)$ for $n\geq 3$
(see \cite{man}).

In this section we will apply our Rips orbit construction to a finitely
generated group $G$ acting by isometries on a quasitree. We will obtain
strong positive results on the structure of $G$ and the way $G$ can act,
usually under some finiteness condition on the action such as having a
locally finite orbit, although Corollary \ref{serc} is completely
general. In the next section we will combine these results with those
in \cite{msw} to turn unbounded actions of $G$ on quasitrees into unbounded
actions on trees, using only the assumption of a locally finite orbit.

A very useful criterion to tell whether a geodesic metric
space $X$ is a quasitree is Manning's bottleneck criterion (see \cite{mn1}).
This states that there is some number $C\geq 0$ (the bottleneck constant)
such that for every geodesic segment $[x,y]$ in $X$ and any point $z$ on
$[x,y]$, any path between $x$ and $y$ must intersect the closed ball
$B(z,C)$. We can use this to examine finitely generated
groups acting on quasitrees by
combining it with the following, which is a straightforward extension of the
point in Section 2 about the connectivity of the Rips orbit graph.
\begin{lem} \label{space}
Let a finitely generated group $G$ with finite generating set $S$ act on
a metric space $(X,d)$ by isometries. For any $x_0\in X$, put
$M=\mbox{max}_{s \in S}d(x_0,s(x_0))$.

Suppose we have partitioned the orbit $G\cdot x_0$ into pairwise disjoint
non empty subsets $A_i$ for $i\in I$. Given $i,j\in I$ with $i\neq j$, set
\[d(A_i,A_j)=\mbox{inf }\{d(a,b)\,|\,a\in A_i, b\in A_j\}.\]
If there is $N>0$ with $d(A_i,A_j)\geq N$ for all $i\neq j$ then $N\leq M$.
\end{lem}
\begin{proof}
Suppose that $d(A_i,A_j)>M$ for all $i\neq j$. If $x_0\in A_1$ say then
$d(s^{\pm 1}(x_0),x_0)\leq M$ and so $s^{\pm 1}(x_0)\in A_1$ too for each
$s\in S$.

We now proceed by induction on the word length of $g\in G$ with respect
to $S$. If all elements of $G$ with word length at most $n$ are in $A_1$
then we can write an element $g$ of word length $n+1$ as $g=g_1s^{\pm 1}$
where $s\in S$ and $g_1(x_0)\in A_1$. This says that 
$d(g(x_0),g_1(x_0))\leq M$ so $g(x_0)\notin A_j$ for $j\neq 1$. Thus
$g(x_0)\in A_1$ and so is the whole orbit of $x_0$.
\end{proof}

We use this as follows.
\begin{prop} \label{quas}
Suppose we have a group $G$ generated by the finite set $S$
which acts by isometries on some quasitree $Q$ with
bottleneck constant $C$. Pick any $x_0\in Q$ and set
$M=\mbox{max}_{s \in S}d(x_0,s(x_0))$.
Then the orbit
$G\cdot x_0$ is $K$-quasiconvex for $K$ satisfying $2K-2C\geq M$.
\end{prop}
\begin{proof}
Take such a $K$ and  
suppose we have $x,y\in G\cdot x_0$ and a geodesic segment $[x,y]$ in $Q$
but there is $z\in [x,y]$ having distance more than $K$ to any point
of $G\cdot x_0$. Removing the ball $B(z,C)$ from $Q$ makes it (path)
disconnected and with $x,y$ in different (path) components. Let us
intersect each of the components of $Q\setminus B(z,C)$ with the
orbit $G\cdot x_0$ (which is contained in $Q\setminus B(z,C)$)
to form the sets $A_i$ for some $i$ in some indexing
set $I$.

On taking arbitrary points $a\in A_i$ and $b\in A_j$ for $i\neq j$ and
a geodesic segment between them, we see that $[a,b]$ has to touch
$B(z,C)$ as these points are in different components on removal of this
closed ball. So suppose $[a,b]$ hits $B(z,C)$ at $p\in [a,b]$. We have
$d(a,z)>K$ and $d(z,p)\leq C$ so $d(a,p)>K-C$ and also $d(b,p)>K-C$.
But we have
\[d(a,p)+d(p,b)=d(a,b)>2(K-C)\geq M.\]
As we have at least two non empty sets $A_i$ and $A_j$ (one containing $x$ and
one $y$), we have a contradiction from Lemma \ref{space}.
\end{proof}   

We now have:

\begin{co} \label{gam}
Suppose a finitely generated group $G$ acts by isometries on a
quasitree $Q$. Then for any $x_0\in Q$
there is $r>0$ such that there is a coarse $G$-equivariant quasi-isometric
embedding from the Rips orbit graph $\Gamma_r(x_0)$ (with its
own path metric) to $Q$.

If $G\cdot x_0$ is a locally finite orbit then $\Gamma_r(x_0)$
has bounded valence.
\end{co}
\begin{proof}
On taking any $x_0\in Q$ we can first apply Proposition \ref{quas} and
then Proposition \ref{os}. We saw in Section 2 that locally finite
orbits give rise to bounded valence Rips orbit graphs.
\end{proof}

We now need to turn this quasi-isometric embedding of $\Gamma_r(x_0)$ in $Q$
into a quasi-isometry to a smaller graph. 
\begin{thm} \label{qism}
Suppose a finitely generated group $G$ acts by isometries on a quasitree $Q$.
Then the Rips orbit graph as obtained
from Corollary \ref{gam} is quasi-isometric to a tree $T$. If this orbit
graph is based at $x_0\in Q$ with $G\cdot x_0$ a locally finite orbit then 
we can take $T$ to have bounded valence.
\end{thm}
\begin{proof}
As $\Gamma_r(x_0)$ quasi-isometrically
embeds in $Q$ which by definition is quasi-isometric
to some simplicial tree $S$, we have that $\Gamma_r(x_0)$ also
quasi-isometrically embeds in $S$ under the map $q$ say. We can assume that
this map sends $\Gamma_r(x_0)$, and in particular its vertices, to vertices
of $S$. We construct $T$ by setting it equal to the convex span in $S$ of
the image of $\Gamma_r(x_0)$ in $T$. Note that for any two points in $T$,
the distance between them in $T$ and in $S$ is the same.

Let $q$ satisfy
\[\frac{1}{C}d_\Gamma(x,y)
-\epsilon\leq d_S(q(x),q(y))\leq Cd_\Gamma(x,y)+\epsilon.\]
We need to show that our map has quasi-dense image in $T$.
Let us take any point $t\in T$: by construction this
means that we have $x,y\in S$ which are the images of vertices $v,w$ say
in $\Gamma_r(x_0)$ and $t$ lies on the geodesic in $S$ between $x$ and $y$.

Removing $t$ splits $S$ into a number (possibly infinite) of path connected
components. First suppose that no point of the image of $q$ lies in the
closed ball $B(t,C+\epsilon)$. Take a geodesic path in $\Gamma_r(x_0)$ from
$v$ to $w$ and consider the images in $S$ of the vertices
$v_0=v,v_1,\ldots ,v_n=w$ making up this path. As $x$ and $y$ lie in
different (path) components $C_1,C_2$ say of $S\setminus\{t\}$, we have $i>0$
such that $q(v_{i-1})\in C_1$ but $q(v_i)$ is not. Now these two points are
outside $B(t,C+\epsilon)$ but the geodesic connecting them must pass
through $t$, thus $d_S(q(v_{i-1}),q(v_i))>2C+2\epsilon$. However
$d_\Gamma(v_{i-1},v_i)=1$ which is a contradiction.

Let us now consider the case where $t$ is as before, lying on the geodesic
$[x,y]$ in $S$, but where $t$ is also a vertex of $T$.
Our argument above will now also show that
for any (path) component of $S\setminus\{t\}$
which intersects $T$, say in the vertex $p\in T$,
there is a vertex $v_t$ of $\Gamma_r(x_0)$
with $q(v_t)$ lying in that component and at distance at most
$C+\epsilon$ from $t$. Suppose this is not true. There must be a vertex $z$
in this path component which is in the image of $q$, since $p$ being
in $T$ implies that it must
lie on a geodesic with both endpoints in the image of $q$,
at least one of which will lie in this path component. We can now take a 
geodesic path in $\Gamma_r(x_0)$ starting from a vertex which is mapped
to $z$ and ending at $x$ which is in a different path component without
loss of generality (because $x$ and $y$ are not in the same path component,
so swap $x$ with $y$ if necessary). As before,
the image of this path at some point has to jump by a distance more than
$C+\epsilon$ in order to get to another path component which is a
contradiction. 

Thus if $G\cdot x_0$ is a locally finite orbit in $Q$, so that $\Gamma_r(x_0)$
has bounded valence, then on taking a base vertex $v_0$ in
\[V:=\{v\in\Gamma_r(x_0)\,|\,d_S(q(v),t)\leq C+\epsilon\},\]
we have for any $v\in V$ that
$d_S(q(v),q(v_0))\leq 2C+2\epsilon$ so $d_\Gamma(v,v_0)\leq 
2C^2+3C\epsilon$. But bounded valence of $\Gamma_r(x_0)$
means that the size of $V$ is
bounded independently of $v_0$ and $t$, thus so is the number of components
of $S\setminus\{t\}$ which intersect $T$. This is an upper bound for the
number of edges in $T$ leaving $t$.
\end{proof}
 
We can now obtain some strong consequences of actions on quasitrees
provided that our group is finitely generated.
\begin{co} \label{act}
Suppose that $G$ is a finitely generated group acting by isometries on an
arbitrary quasitree $Q$. Then there is
a quasitree $Q_0$ which is a graph equipped with its path metric 
where $G$ acts coboundedly by graph automorphisms and a coarse $G$-equivariant
quasi-isometric embedding from $Q_0$ to $Q$.
\end{co}
\begin{proof}
Form the Rips orbit graph $Q_0=\Gamma_r(x_0)$ by taking any $x_0\in Q$
and choosing $r>0$ as in Corollary \ref{gam}. Then $G$ acts on it
by graph automorphisms and coboundedly.
By Theorem \ref{qism} this graph is quasi-isometric to a tree.
\end{proof}

Consequently if the original action on $Q$ possesses the properties we
have discussed, such as having unbounded orbits, or is of general type,
or metrically proper, then so does the resulting action of $G$ on $Q_0$.
We can take this further:
\begin{co} \label{propq}
If $G$ is a finitely generated group acting by isometries
and metrically properly on an arbitrary quasitree $Q$ then $G$ is
a virtually free group. 
\end{co}
\begin{proof}
Apply Corollary \ref{act} so that $G$ acts on $Q_0$ by isometries,
metrically properly and coboundedly. By Svarc - Milnor
this implies that our group $G$ is quasi-isometric to $Q_0$. But
all orbits are locally finite so $Q_0$ has bounded valence and by
Theorem \ref{qism} so does $T$, which is also quasi-isometric to $G$.
As $T$ and $G$ are proper hyperbolic spaces, their boundaries must
be homeomorphic and so $G$ is a virtually free group (see \cite{ikpbn}
Proposition 2.20 (b) and Theorem 8.1).
\end{proof}
Note that by the comment at the end of Section 2, we actually have that
the orbit maps from $G$ to $Q_0$ and to $Q$ must be a quasi-isometry
and a quasi-isometric embedding respectively.

\begin{co} \label{bndd}
Suppose that $G$ is a finitely generated group acting by isometries and with
unbounded orbits on a quasitree $Q$ and there is a locally finite orbit.
Then $G$ acts coboundedly, by graph automorphisms and with unbounded orbits
on a bounded valence quasitree $Q_0$ which is quasi-isometric to a
bounded valence tree.
\end{co} 
\begin{proof}
As before but now by taking a basepoint $x_0\in Q$ where $G\cdot x_0$
is a locally finite orbit, the corresponding
Rips graph has bounded valence and it is quasi-isometric to a tree which
also has bounded valence by Theorem \ref{qism}.
\end{proof}
\begin{co} If $G$ is a finitely generated group acting by graph automorphisms
on a locally finite graph $Q$ which is a quasitree and the stabilisers are
finite then $G$ is a virtually free group. In particular if $G$ acts
freely then it is a virtually free group.
\end{co}
\begin{proof} The finite stabilisers condition when applied to a group
acting by graph automorphisms means that $G$ acts topologically properly
on $Q$. But as $Q$ is locally finite, this will be
a metrically proper action so Corollary \ref{propq} applies.
\end{proof}

We finish this section by looking at isometric actions of a group on a
quasitree where every element is elliptic. Of course this could happen
because the action has bounded orbits. But it could also be a parabolic
action even in the case of a tree, indeed even in the case of a metrically
proper action on a bounded valence tree, as for the coset
construction. However these
groups are never finitely generated and we mentioned in the examples at the
start of this section that if a finitely generated group $G$ acts on
a tree with every element elliptic, then the action has bounded orbits
by \cite{ser} I.6.5. This result of Serre's does not hold for
actions on quasitrees. Nevertheless we can still use
our arguments to show that finitely generated groups have no parabolic actions
on quasitrees.
\begin{co} \label{serc}
If $G$ is a finitely generated group acting by isometries
on an arbitrary quasitree $Q$ then the action
cannot be parabolic. In particular if   
every element of $G$ is elliptic then
$G$ has bounded orbits.
\end{co}
\begin{proof}
We apply Corollary \ref{gam} to find $r>0$ which gives us our coarse
$G$-equivariant quasi-isometric embedding from the Rips
orbit graph $\Gamma_r(X_0)$ to $Q$. Thus the hyperbolic type of the two
actions are the same, but the action of $G$ on the Rips orbit graph
is cobounded so it cannot be parabolic, thus nor is the action of
$G$ on $Q$.
\end{proof}
A result in \cite{man} Section 3 is that any element $g$ acting by an isometry
on a quasi-tree is either elliptic or loxodromic. By applying the above
to $\langle g\rangle$ we recover this result.

\section{Unbounded orbits on quasitrees}

We have mentioned many examples of groups $G$ acting on quasitrees in various
ways but where there exists no action of $G$ of the same type on a tree.
However we saw in the last section that if $G$ is finitely generated
then there are restrictions on how $G$ can act on
a quasitree, especially if there is some sort of finiteness condition
on the action. In particular, Corollary
\ref{propq} tells us that under the finite generation condition, $G$
acting metrically properly on a quasitree implies that $G$ is virtually
free and so acts metrically properly on a tree.

In \cite{msw} a strong result was obtained about arbitrary 
cobounded quasi-actions of an 
arbitrary group on trees satisfying the finiteness property of having
bounded valence:
\begin{thm} \label{mswt} (\cite{msw} Theorem 1)\\
If $G\times T\rightarrow T$ is a cobounded quasi-action of a group $G$
on a bounded valence bushy tree $T$, then there is a bounded valence, bushy
tree $T'$, an isometric action $G\times T'\rightarrow T'$, and a
coarse $G$-equivariant quasi-isometry from the action of $G$ on $T'$ to
the quasi-action of $G$ on $T$.
\end{thm}
Here a {\it bushy} tree $T$ is one where there exists $b\geq 0$ (the
bushiness
constant) such that for any point $x\in T$,
$T\setminus B(x,b)$ has at least three unbounded
components.
Meanwhile a $(C,\epsilon)$ {\it quasi-action} of a group $G$ on
a metric space $X$ is a map $\alpha:G\times X\rightarrow X$ such that
for each $g\in G$ the map $A_g$ given by $A_g(x)=\alpha(g,x)$ is a
$(C,\epsilon)$ quasi-isometry and for each $g,h\in G$ the supremum
over $x\in X$ of the distance in $X$ between $\alpha(g,\alpha(h,x))$
and $\alpha(gh,x)$ is at most $C$. The key point for us here is that
if $G$ acts by isometries on $X$ and $q:X\rightarrow Y$ is a quasi-isometry
with a quasi-inverse $r:Y\rightarrow X$ then 
$\alpha(g,y):=q(g(Q(y)))$ is a quasi-action of $G$ on $Y$ and then $q$
and $r$ are both coarse $G$-equivariant.

We also say that a quasi-action $\alpha$ is {\it cobounded} 
if there is $R\geq 0$
such that for each $x_0\in X$, any $x$ in $X$ is at most a distance $R$
to some point in the quasi-orbit $\{\alpha(g,x_0)\,|\,g \in G\}$. This
notion is preserved by coarse $G$-equivariant quasi-isometries
and reverts to its usual meaning if the quasi-action of $G$
on $X$ is a genuine isometric action.

Thus looking at our Corollary \ref{bndd}, for a finitely generated group
$G$ acting with unbounded orbits on a quasitree $Q$ and having a locally
finite orbit, we obtain from this
a cobounded action of $G$ on a quasitree $Q_0$ which also has unbounded
orbits and where $Q_0$ is quasi-isometric to a bounded valence tree
$T_0$. We can then turn our action of $G$
on $Q_0$ into a quasi-action of $G$ on $T_0$ which is also cobounded.
If we knew that $T_0$ was a bushy tree then Theorem \ref{mswt} says
that $G$ acts isometrically and coboundedly
on a bounded valence tree $T_0'$. This action will have unbounded orbits
because $T_0'$ is bushy, so is not a finite tree. (Indeed this argument 
works on replacing unbounded orbits with any notion that can be defined
for quasi-actions, is preserved under a coarse $G$-equivariant quasi-isometry
between two quasi-actions, and is preserved when taking the Rips orbit
graph.)

Therefore in this section we will apply Theorem \ref{mswt} by starting
with an appropriate action of a finitely generated group $G$ on a
quasitree and finding out when the corresponding tree $T_0$ is bushy,
as well as seeing what conclusions we can draw when it is not.
We first consider when we have a bushy tree. 
\begin{prop} \label{mswp}
Suppose that a group $G$ acts by graph automorphisms and vertex transitively
on a bounded valence graph $Q_0$ which
is quasi-isometric to a bounded valence tree $T$. If $Q_0$ has at least 
three ends then $T$ is a bushy tree.
\end{prop}
\begin{proof}
As $Q_0$ and $T$ are both locally finite, all the various notions
of ends coincide in this case. For a tree, the set of ends and
the boundary $\partial T$ are the same and this is true for $Q_0$ too
because of the quasi-isometry, so that $T$ and $Q_0$ have the same
number of ends.

As $T$ has at least three ends, 
there
is some vertex $v_0\in T$ such that $T\setminus\{v_0\}$ has
at least three unbounded components (and possibly some bounded ones too).
Take
three sequences $(u_n),(v_n),(w_n)$ of adjacent vertices in $T$, all distinct
apart from $u_0=v_0=w_0$, such that each sequence tends to a different
end of $T$. Note that we have
$d(u_m,v_n)=d(v_m,w_n)=d(w_m,u_n)=m+n$ for any integers $m,n\geq 0$,
where $d$ denotes the usual distance in $T$.

Now we have a quasi-isometry $q:Q_0\rightarrow T$ so
there is $C>0$, not depending on any choice made so far and which we will
take to be an integer, such that
we have a vertex $p_0\in Q_0$ with its image $q(p_0)$ within distance $C$
of $v_0$ in $T$.

If $T$ were not bushy then for any integer $M> C$
we would have a vertex $x_M$ of
$T$ such that $T\setminus B(x_M,M)$ has at most two unbounded components.
But there will also exist a vertex $p_M\in Q_0$ such that
$q(p_M)$ lands in the ball $B(x_M,C)$. As the action of $G$ on
$Q_0$ is isometric and vertex transitive, we have some isometry $g_m\in G$
(or at least its image in $Aut(Q_0)$ if this action is not faithful)
with $g_M(p_0)=p_M$. But by taking a quasi-inverse $r$ of $q$ with
$p_0=r(v_0)$, we can
``quasi-conjugate'' $g_M$ into a quasi-isometry $h_M:=qg_Mr$ of $T$ to itself.
Moreover $h_M$ is a $(K,\epsilon)$ quasi-isometry for $K\geq 1$ and
$\epsilon\geq 0$
independent of $M$ because $g_M$ is an isometry. We will henceforth
assume that $K$ and $\epsilon$ are integers.

Let $(a_n), (b_n), (c_n)$ be the three $(K,\epsilon)$ quasi-geodesics which
are the images under $h_M$ of the geodesics $(u_n),(v_n),(w_n)$
respectively.
These all start at $h_M(v_0)=q(p_M):=z_0$
    which lies inside $B(x_M,C)$. Now if $s_1$ and $s_2$ are the unique
    nearest points of the unbounded components of
    $T\setminus int(B(x_M,M))$
    from $x_M$ (with possibly $s_1=s_2$ if these two components are the
    same) then each of our sequences $(a_n),(b_n),(c_n)$ must end up
    in one of these unbounded components. Assuming without
    loss of generality that $(a_n)$ and $(b_n)$
    end up in the same one, say the one with endpoint $s_1$,
     let us take some large $N$ with both the vertices
     $a_N$ and $b_N$ in this component. Consider the two geodesics
$\gamma,\delta$
     from
     $z_0=a_0$ to $a_N$ and from $z_0=b_0$ to $b_N$.
Let us label
the vertices of these geodesics $\gamma_0=z_0,\gamma_1,
\ldots ,\gamma_A=a_N$
     and $\delta_0=z_0,\delta_1,\ldots ,\delta_B=b_N$
     where $A=d(z_0,a_N)$ and $B=d(z_0,b_N)$. Note that
     they are the
     same geodesic $\eta$ from $z_0$ to the point $s_1$, so
     $\gamma_i=\delta_i$ for $0\leq i\leq L$, where
     where $L=d(z_0,s_1)$ which is at least $M-C$.
     Now $a_0,a_1,\ldots ,a_N$ and $b_0,b_1,\ldots ,b_N$ are both
     $(K,\epsilon)$ quasi-geodesic segments
     (namely a $(K,\epsilon)$ quasi-isometric embedding from $\{0,1,\ldots ,N\}$
     to $T$), so a version of the Morse lemma for trees tells us that
     there is some number $R\geq 0$ depending only on $K$ and $\epsilon$
     such that every vertex $a_i$ on the first quasi-isometric segment
     is at most $R$ from some vertex lying on
     $\gamma$, and the same for the
     second segment and $\delta$.

     Let us pick $M$ and hence $L$ very large and consider any $i$ with
     $i\geq 0$ and $i+K+\epsilon\leq L-1$. Consider the geodesic segment
     from $\gamma_0$ to $\gamma_L$
     which is common to both $\gamma$ and $\delta$ and
     union it with the $K+\epsilon$ rooted trees in $T$ which have roots
     $\gamma_{i+1},\ldots ,\gamma_{i+K+\epsilon}$ but not including any edges of this
geodesic segment. Removal of this whole subtree results in two components
     separated by a distance greater than $K+\epsilon$, but as we have
     $d(a_i,a_{i+1})\leq K+\epsilon$ with $a_0$ in one component and
     $a_N$ in another, we must have some $a_j$ lying in this subtree.
     In exactly the same way, there is some $b_k$ lying in there too.
     By application of the Morse lemma, we obtain
     $1\leq s,t\leq K+\epsilon$ such that
     \[d(a_j,\gamma_{i+s}) \leq R\mbox{ and } d(b_k,\gamma_{i+t}) \leq R\]
     but
     \[ d(\gamma_{i+s},\gamma_{i+t})\leq K+\epsilon-1\mbox{ thus }
     d(a_j,b_k)\leq 2R+K+\epsilon-1.\]

     But as the geodesics from $\gamma_0$ to $a_j$ or to $b_k$ each pass
     through the vertex $\gamma_i$, we have
     \[Kj+\epsilon\geq d(a_j,a_0)\geq d(\gamma_i,\gamma_0)=i\]
     and also $Kk+\epsilon\geq i$. Thus $j+k\geq 2(i-\epsilon)/K$.
     Now originally we had $d(u_j,v_k)\geq j+k$ and the sequences $(a_n)$,
     $(b_n)$ we obtained by applying the $(K,\epsilon)$ quasi-isometry
     $h_M$ to $(u_n),(v_n)$. Thus
     \[2R+K+\epsilon-1\geq d(a_j,b_k)\geq (j+k)/K-\epsilon
     \geq 2(i-\epsilon)/K^2-\epsilon.\]
     But here the left hand side only depends on $K,\epsilon$ which are
     constant throughout the proof. Thus by letting $M$ tend to
     infinity we have
     $L$ at least $M-C$, with $C$ also constant throughout so $L$ tends
     to infinity too. We can then take $i=L-K-\epsilon-1$ to obtain a valid
     $i$ in the right hand side of the equation above, which must therefore
     tend to infinity too with $M$, giving us a contradiction.
     \end{proof}

By results in \cite{ser} we know that a finitely generated group has
an unbounded action on a tree (equivalently has an action without a
global fixed point) if and only if it splits as an amalgamated
free product or HNN extension.
We can now find out which finitely generated groups act with unbounded orbits
on quasitrees under the assumption that the action has a locally finite
orbit.
\begin{thm} \label{summ}
Let $G$ be any finitely generated group having an unbounded action by
isometries on some quasitree $Q$. Suppose that there is a locally
finite orbit.
Then there is a bounded valence tree $T'$ on which $G$ acts coboundedly
by automorphisms and a coarse $G$-equivariant quasi-isometric
embedding from $T'$ to $Q$.
\end{thm}
\begin{proof}
We consider the hyperbolic type of the action of $G$
on $Q$. By hypothesis this is not bounded and by
Corollary \ref{serc} it cannot be parabolic.
If this action is lineal with limit points $\partial_GQ=\{x^+,x^-\}$
then by Corollary \ref{cgqie} we can take $T'$ to be the simplicial
line $L$.

In the two cases where the action is either quasi-parabolic or of general
type, we apply Corollary \ref{gam} to our action and obtain the
graph $\Gamma$ on which $G$ acts by automorphisms and vertex transitively,
along with a coarse $G$-equivariant
quasi-isometric embedding from $\Gamma$ to $Q$. But $\Gamma$ has
bounded valence because of the locally finite orbit and moreover
by Theorem \ref{qism}
it is quasi-isometric to a bounded valence tree $T$.

Now the hyperbolic type of the action is preserved when passing from
$Q$ to $\Gamma$. But a quasi-parabolic or general type action on $\Gamma$
implies that the limit set $\partial_G\Gamma$ has at least three points,
so the boundary $\partial\Gamma$ has at least three points too and hence at
least three ends because it is a quasitree.
Hence by Proposition \ref{mswp} $T$ is a bushy tree.
Thus application of Theorem \ref{mswt} gets us a
bounded valence bushy tree $T'$, an action of $G$ on $T'$ which is isometric
and thus here by automorphisms, 
and a coarse $G$-equivariant quasi-isometry from the action of $G$ on $T'$
to the quasi-action of $G$ on $T$. But this quasi-action was obtained
by taking the action of $G$ on $\Gamma$ and quasi-conjugating, so we also
have a coarse $G$-equivariant quasi-isometry from $T$ to $\Gamma$
and we started with a coarse $G$-equivariant quasi-isometric embedding
from $\Gamma$ to $Q$. Thus on
composing these maps from $T'$ to $T$ to $\Gamma$ to $Q$, we end up with a
coarse $G$-equivariant quasi-isometric embedding from the action of $G$ on
$T$ to the original action of $G$ on $Q$.
\end{proof}

So in particular if the original action of $G$ on $Q$ is metrically
proper as in Corollary \ref{propq} then, as before, we see
that $G$ acts geometrically by
automorphisms on $T'$ with the orbit maps from $G$ to $T'$ and $\Gamma$
both quasi-isometries and the orbit map from $G$ to $Q$ a
quasi-isometric embedding.

As any action by automorphisms on a locally finite graph has all orbits
locally finite, we immediately obtain:
\begin{co} \label{qttt}
  Suppose that $G$ is a finitely generated group acting by
automorphisms
and with unbounded orbits on a locally finite graph which is a quasitree.
Then $G$ acts by automorphisms
and without a global fixed point on some simplicial tree with bounded valence.
Thus $G$ has the structure of a fundamental group of a (non trivial) finite
graph of groups where all edge groups have finite index in all of the
respective vertex groups.
\end{co}
\begin{proof}
The first part is application of Theorem \ref{summ} to obtain our action of $G$
on a bounded valence tree $T'$. If this had bounded orbits, which for a
tree is equivalent to having a global fixed point, then so would the
original action.

The statement on the structure of $G$ is a standard application of Serre's
results in \cite{ser}, taking care to subdivide the edges of $T'$ (but
leaving $G$ unchanged) if any elements of $G$ invert an edge.
\end{proof}

Having mentioned in the last section that if a finitely generated group
acts on a quasitree with various elements acting elliptically then this
is not enough to force a bounded orbit, unlike the case of a tree, we
can at least get this to work on locally finite quasitrees.
\begin{co}
  Suppose we have a theorem that says: let $G$ be a finitely generated
  group acting on a tree where certain elements are elliptic, then $G$
  has a global fixed point (such as \cite{ser} I.6.5). Then for any
  action of $G$ on a locally finite graph $Q$ which is a quasitree, if those
  same elements are elliptic on $Q$ then the action of $G$ on $Q$ has
  bounded orbits.
\end{co}
\begin{proof}
  By applying Theorem \ref{summ} to $G$ acting on $Q$, we have our action
  of $G$ on the tree $T'$. But this is of the same hyperbolic type
  as our original action and with elements acting in the same way,
  so the elements which were elliptic on $Q$ are still elliptic on $T'$
  and thus our theorem in the statement applies. So the action of $G$
  on both $T'$ and $Q$ has bounded orbits.
\end{proof}

We finish this section by giving a counterexample to Theorem \ref{summ}
when the group is not finitely generated, even when the action is
metrically proper.\\
\hfill\\
{\bf Example}: Consider the restricted direct product
$G=C_2\times C_3\times C_5\times C_7\times \ldots$. By regarding $G$ as a
direct union of the increasing sequence of subgroups
$C_2,C_2,\times C_3,C_2\times C_3\times C_5,\ldots $ we can use the coset
construction to obtain a parabolic action of $G$ on a locally
finite tree but with unbounded valence. As all stabilisers are finite,
this is a proper action.

Now suppose $T$ is any tree with
valence bounded by $N$ and consider an element $g$ acting on $T$ with
$g$ having order coprime to the numbers 
$2,3,4,\ldots ,N$. Then $g$ must fix a vertex $v_0$ as
$T$ is a tree but it must then also fix the vertices adjacent to $v_0$
by Orbit - Stabiliser and so on, thus $g$ acts as the identity. This means
that if $p^+$ is the smallest prime which is greater than $N$
and $p^-$ the largest prime less than $N$
then any element
in the subgroup $C_{p^+}\times\ldots$ acts trivially on $T$. Thus the quotient
of $G$ by this subgroup which is in the kernel of the action leaves us
only the finite group $H=C_2\times C_3\times\ldots \times C_{p^-}$
and so any action of $G$ will have bounded orbits.

\section{Acylindrical actions}

We say that an isometric action of a group $G$ on a hyperbolic space $X$ 
is {\it acylindrical}
if given $\epsilon\geq 0$ there is $N$ and $R$ such that for any two points
$x,y\in X$ at least a distance $R$ apart, there are no more than $N$
elements $g\in G$ moving both $x$ and $y$ at most $\epsilon$. We note
some basic consequences of this definition: if $X$ is a bounded metric
space then it is automatically satisfied, as well as if $G$ is finite.
In the unbounded case the action need not be faithful
but the kernel of the action must be finite. However an 
action on an unbounded metric space
with bounded orbits may or may not be acylindrical. We also note here that
on being given a coarse $G$-equivariant quasi-isometric embedding from
$Y$ to $X$, acylindricity of the action on $X$ passes to $Y$ but
this need not be true when going from $Y$ to $X$.

The definition of a {\it acylindrically hyperbolic} group is one that admits
an isometric action on a hyperbolic space which is acylindrical but not
an elementary action. This has notable consequences, particularly that the
group is SQ-universal. A summary of results and properties of acylindrical
groups is provided in \cite{osac}. In particular we note that there are
strong restrictions on elementary acylindrical actions: they cannot be
parabolic or quasi-parabolic. Therefore there are no parabolic elements
in any acylindrical action, because acting acylindrically is clearly
preserved by subgroups. Moreover for a group action which is acylindrical
and non elementary, restricting the action to any subgroup
$H$ is also acylindrical and
non elementary, or it has bounded orbits, or it has two limit points. In
the latter case $H$ will be a virtually cyclic group
(moreover all virtually cyclic
groups with a loxodromic element will act acylindrically, see the
next section).

Many groups, such as mapping class groups and (outer) automorphisms of
free groups (ignoring finitely many small examples), are known to be
acylindrically hyperbolic, as is of course any word hyperbolic group $G$
by utilising the natural action on its own Cayley graph $C$ which is
an acylindrical action. On considering
the action by restriction of any subgroup $H$ of $G$, we have that this action
will only have bounded orbits if $H$ is finite. Thus if $H$ is not virtually
cyclic then it too will be acylindrically hyperbolic; moreover it has a
non elementary acylindrical action not just on an arbitrary hyperbolic space
but on a proper hyperbolic space.
However for mapping class groups, one uses
the curve complex to obtain an acylindrical action and this is a graph
which is not locally finite and so not a proper metric space. Therefore
it is reasonable to ask which groups have a non elementary acylindrical action on a
{\it proper} hyperbolic space and this turns out to be much more
restrictive than the class of acylindrically hyperbolic groups.
Indeed, if we look at the case where the hyperbolic space
is a locally finite (but unbounded) 
graph then, as noted in \cite{balas}, 
an acylindrical action of a group $G$ cannot have
stabilisers of infinite order. To see this, note that if a vertex $v_0\in V$
is fixed by the infinite subgroup $G_0$ of $G$ then for each vertex $v$ of
distance 1 away from $v$, we have a finite index subgroup of $G_0$ fixing
$v$ too. Thus any two vertices have a common stabiliser that is infinite
so the action is not acylindrical (unless it is a finite graph). This
means that $G$ acts on the graph metrically properly, which we have
seen gives strong restrictions on $G$.

The next theorem is an extension of this to any proper hyperbolic space.
In this section we will write proper action for what we have called a
metrically proper action, as it coincides with other notions of a proper
action when the metric space is proper.
\begin{thm} \label{acyl}
Suppose that an arbitrary group $G$ acts acylindrically on any proper
metric space $X$. Then either $X$ is bounded or $G$ acts properly.
\end{thm}
\begin{proof}
Suppose that the action is not proper, so that there is $x_0\in X$
and $r>0$ with the closed ball $B(x_0,r)$ containing $g(x_0)$ for
infinitely many $g$. Let $S=\{g\,|\,g(x_0)\in B(x_0,r)\}$ be this infinite
set.
Now take any $R>0$ and any $y$ with $d(x_0,y)\geq R$
(if there is no $y$ then $X$ is bounded). Note that for any $g\in S$ we have
\[d(g(y),x_0)\leq d(g(y),g(x_0))+d(g(x_0),x_0)\leq d(x_0,y)+r:=k\]
where $k$ is independent of $g$. 
So as $g(y)$ lies in the compact set $B(x_0,k)$, we can take a sequence of
distinct elements $g_n\in S$ such that $g_n(y)$ tends to some point $z$.
Then
\[d(g_n^{-1}g_{n+1}(y),y)\leq d(g_{n+1}(y),y)+d(g_n(y),y)\]
and so $g_n^{-1}g_{n+1}(y)$ tends to $y$. But
\[d(g_n^{-1}g_{n+1}(x_0),x_0)\leq d(g_{n+1}(x_0),x_0)+d(g_n(x_0),x_0)\leq r+r\]
so that for sufficiently large $n$ the elements $g_n^{-1}g_{n+1}$ of $G$ move
both $x_0$ and $y$ at most $2r$. 
If $g_n^{-1}g_{n+1}$ made up only
finitely many elements of $G$ then $g_n^{-1}g_{n+1}(y)$ is one of only
finitely many points, so is eventually constant as it tends to $y$.
Without loss of generality it is constant, so $g_n(y)=g_1(y)=y$. But
$g_1^{-1}g_n$ is certainly a sequence of distinct elements which fix $y$
and it also moves $x_0$ by at most $2r$ so the action is not acylindrical
in this case too.
\end{proof} 

Therefore we immediately have by this result and the fact
that acylindrical actions are never parabolic or quasi-parabolic:
\begin{co}
Suppose that an arbitrary group $G$ acts acylindrically on any proper
hyperbolic space $X$ which is not bounded. Then for any subgroup $H$
of $G$ we have exactly one of the following:\\
(i) $H$ acts on $X$ with bounded orbits, which happens if and only if
$H$ is finite\\
(ii) $H$ acts on $X$ with two limit points, which happens if and only
if $H$ is virtually cyclic and infinite.\\
(iii) The action of $H$ on $X$ is acylindrical, proper and non elementary, so
that $H$ is an acylindrically hyperbolic group.
\end{co}
This statement differs from the equivalent result for arbitrary hyperbolic
spaces in (i), where $H$ must be finite here, and also in (ii) where
$H$ being virtually cyclic here implies that it must act with unbounded orbits.

We thus immediately obtain the following alternative.
\begin{co}
Let $G$ be any group acting acylindrically on a proper hyperbolic space
which is unbounded. Then $G$ is either acylindrically hyperbolic or
virtually cyclic and the same applies for any subgroup of $G$. In
particular $G$ does not contain $\Z\times\Z$.
\end{co}
So for instance, the initial examples of groups which are not word hyperbolic
but are relatively hyperbolic were the fundamental groups of
finite volume, but not cocompact, complete hyperbolic 3-manifolds. These
are certainly acylindrically hyperbolic groups but this will never be
witnessed by any action on a proper hyperbolic space.

If we stick with proper hyperbolic metric spaces $X$ but now ask conversely
whether a group action on $X$ which is metrically proper will also be
acylindrical,
or whether this group $G$ will act
acylindrically on some other hyperbolic space,
we cannot of course say anything in the case of a parabolic action as then
$G$ could be any countable group. However if this action has a loxodromic
element
$g\in G$ then it is noted in \cite{osac} Section 8(c) that $g$ will be a
WPD element for this action of $G$ and so $G$ is either virtually
cyclic or an acylindrically hyperbolic group. But as we have already seen
for cusped hyperbolic 3-manifolds, $G$ may not have any non elementary
acylindrical action on any proper hyperbolic space.

In \cite{balas} it was shown that every acylindrically hyperbolic group $G$
admits an acylindrical and cobounded non elementary
action on some graph which is a quasitree. It was pointed out that if
we further assume that the graph is locally finite then $G$ must be
(finitely generated and) virtually free.
The final result in this section generalises this statement from quasitrees
which are locally finite graphs to quasitrees which are proper metric spaces
and it does not need the action to be cobounded. It thus displays a big
difference from the situation where the quasitree need not be proper,
at least in the case of finitely generated groups.
\begin{co}
Let $G$ be a finitely generated group having an acylindrical action
on a proper hyperbolic metric space which is unbounded and also a quasitree.
Then $G$ is a virtually free group.
\end{co}
\begin{proof}
We apply Theorem \ref{acyl} and then Corollary \ref{propq}.
\end{proof}

\section{Bounded valence hyperbolic graphs}

We have already mentioned the result in \cite{man} Section 3 that
any isometry of a quasitree is either elliptic or loxodromic.
Another important class of  hyperbolic spaces where any isometry
is either elliptic or loxodromic is that consisting of
hyperbolic graphs which are not just
locally finite but of bounded valence.
To see this holds from the relevant results in the literature, we will
recall yet another notion of proper for a group $G$ acting isometrically on
a metric space $X$: we say the action is {\bf uniformly
proper} if for all $r\geq 0$ there is some $N_r\in\N$ such that for any $x\in X$
the set $\{g\in G:d(g(x),x)\leq r\}$ has size at most $N_r$. This is clearly
preserved by subgroups and is
a strengthening both of an acylindrical and of a metrically proper action.
Note first
that if $X$ is an arbitrary metric space $(X,d)$ and $g$ acts on $X$
loxodromically
then this action is uniformly
proper because $|n|\tau(g)\leq d(g^n(x),x)$
and $\tau(g)$ is independent of $x$. Moreover if a group $G$ acts on $X$
and the restriction of the action to some finite index subgroup is
uniformly proper then the action of $G$ is uniformly proper too by
the pigeon hole principle. Thus any virtually cyclic group will also
act uniformly properly on any metric space provided it contains
a loxodromic element.

We now quote \cite{osac} Lemma 13.4 (credited to Ivanov-Olshanskii)
which states that a uniformly proper action of a group on a hyperbolic
space $X$ cannot be a parabolic action. In particular an isometry $g$ of $X$
cannot be a parabolic element if $\langle g\rangle$ acts uniformly
properly.

Here we specialise to the case of a connected graph $\Gamma$ equipped
with the path metric and with valence bounded by $c$. Given an
automorphism $g$ of $\Gamma$, the action of $\langle g\rangle$
must be uniformly proper on $\Gamma$ if it is infinite
and acts freely. This is
because for any vertex
$v\in V(\Gamma)$ there are only finitely many vertices at distance at
most $R$ from $v$, with this number independent of the actual vertex
chosen. So if further we have that $\Gamma$ is a hyperbolic graph then
in this case $g$ is loxodromic exactly when
$\langle g\rangle$ is infinite and acts freely, and otherwise $g$ will be
elliptic. Note the contrast here with bounded valence
graphs which are not hyperbolic; for instance a distorted cyclic subgroup
$\langle g\rangle$ in a finitely generated group $G$ will act uniformly
properly on a Cayley graph of $G$ but $g$ will be a parabolic
element.

We note here that very recently it was shown in \cite{dmht} using Helly
graphs that a finitely generated group $G$ acting on any bounded valence
hyperbolic graph with every element elliptic
must act with bounded orbits. Thus $G$ cannot act parabolically on such
a space, because the lack of parabolic elements here means that every
element would be elliptic in such an action. This gives us for bounded
valence hyperbolic graphs the exact equivalent of our Corollary
\ref{serc}.

We now look at the structure of groups acting properly or uniformly
properly
on some bounded valence hyperbolic graph. It makes sense to divide this into
the two cases where stabilisers are finite and where there is a uniform bound
to the size of the stabilisers. Note that for a bounded valence (or even
a locally finite) graph, if one point has a finite stabiliser under some
isometric action of a group then all points do.

\begin{prop} \label{bnded}
If a group $G$ acts by automorphisms
on a bounded valence hyperbolic graph
$\Gamma$ with some point having a finite stabiliser than exactly one of
the following occurs:\\
\hfill\\
(i) $G$ is finite\\
(ii) $G$ is infinite but virtually cyclic and contains a loxodromic element.\\
(iii) $G$ is not finitely generated, is locally finite and acts parabolically
with all elements elliptic.\\
(iv) The action of $G$ on $\Gamma$ is general.
\end{prop}
\begin{proof} This will be a metrically proper action, so if it has bounded
  orbits then the group is finite and if it is lineal then the group is
  virtually cyclic. Moreover there are no quasi-parabolic actions which
  are metrically proper.
  
For
a parabolic action, we divide into the cases when $G$ is finitely generated
or not. Proposition 7.1 of \cite{ham} states that if a finitely generated
group $G$ acts metrically properly on a proper hyperbolic space of bounded
geometry, as is the case here, and this action is elementary then
$G$ must be virtually nilpotent. But in a parabolic action,
no infinite order element acts loxodromically whereas in a metrically
proper action, no infinite order element acts elliptically and on a bounded
valence hyperbolic graph there are no parabolics. Thus
$G$ is a finitely generated virtually
nilpotent torsion group and so is finite.

If $G$ acts parabolically but
is not finitely generated then again it cannot contain infinite
order elements and it cannot contain any finitely generated subgroups $H$
that are
infinite by the above argument which also applies to $H$.
\end{proof}

Thus there are strong restrictions on groups
having a proper action by automorphisms
on a bounded valence hyperbolic graph, especially in terms of their subgroup
structure because Proposition \ref{bnded} also applies to their subgroups.

Note that
every countable group $C$ has a proper action on a locally finite hyperbolic
graph (we embed $C$ in a finitely
generated group $G$ and then use the action on the hyperbolic cone) and also
on a bounded valence graph (a Cayley graph of $G$). We mention now
that there are even more extreme
examples.
\begin{co} There exists an infinite finitely generated group whose only
  isometric action on any bounded valence hyperbolic graph, or any locally
  finite graph which is a quasitree, is the trivial
action, where all elements act as the identity.  
\end{co}
\begin{proof} We take a Tarski monster, that is an infinite group $G$ which
is finitely generated but with all proper subgroups finite. Consider any
action of $G$ on a bounded valence hyperbolic graph $\Gamma$ and take
any point $p$ in $\Gamma$. If the stabiliser of $p$ is finite then we
can apply Corollary \ref{bnded}, but $G$ is not virtually cyclic and it
does not contain $F_2$, so cannot have a general action on any hyperbolic
space. So the stabiliser of $p$ is the whole group, but this applies to
all points $p$ in $\Gamma$. As for locally finite quasitrees,
we have by Corollary \ref{qttt} that any action by $G$ would have bounded
orbits, so each point has a finite orbit and thus a finite index stabiliser
which is the whole group.
\end{proof}

We now examine the case where a group acts on some bounded valence
hyperbolic graph with an overall bound on the size of any
stabiliser, as this is the equivalent of being uniformly proper for
such actions.
\begin{prop} \label{upac}
Suppose the group $G$ acts by automorphisms on some unbounded
hyperbolic graph $\Gamma$ whose valence is bounded above by $c$.
The following are equivalent.\\
\hfill\\
(i) All stabilisers are finite and there is a uniform bound on the size
of any stabiliser.\\
(ii) $G$ acts uniformly properly on $\Gamma$.\\
(iii) $G$ acts acylindrically on $\Gamma$.
\end{prop}
\begin{proof}
If (i) holds with $B$ an upper bound for the size of any stabiliser
than let us take any $r\geq 0$ and set $m$ to be the integer below $r$. 
For any vertex $v\in\Gamma$ we
consider the elements $g\in G$ such that $d_\Gamma(g(v),v)\leq r$.
There are at most $1+c+c^2+\ldots +c^m$ vertices within distance $r$
of $v$ and at most $B$ group elements $g\in G$ sending $v$ to the
same vertex, so we have a bound on the number of such group elements
which is independent of the chosen vertex $v$.

Then (ii) implies (iii) is clear. Now suppose (i) does not hold, so that
for any integer $M$ there is a vertex $v_M$ in $\Gamma$ which is stabilised
by at least $M$ elements. But then by Orbit - Stabiliser we have that each
neighbouring vertex of $V_M$ is stabilised by at least $M/c$ elements
that also stabilise $v_M$, each
vertex distance 2 from $V_M$ is stabilised by at least $M/c^2$ elements
that also stabilise $v_M$
and so on. Thus on being given $N,R>0$ we set $s$ to be the integer above
$R$ and find a vertex $v_0$ fixed by at least $Nc^s$ elements. Then any vertex
which has distance $s$ from $v_0$ is fixed by at least $N$ elements which
also fix $v_0$, thus the action of $G$ is not acylindrical.
\end{proof}

We finish this section by noting some examples in the literature of groups
possessing a uniformly proper action on some hyperbolic space. This class
is considered in \cite{coos} where it is denoted by
$\cal P$. They also consider the subclass ${\cal P}_0$ of all groups with
a uniformly proper action on a bounded valence hyperbolic graph, which is
the same as the groups in Proposition \ref{upac}. This in turn contains
the class $\cal S$ of all subgroups of word hyperbolic groups.
They did not know whether all three classes of groups are the same
but gave an example of a group in ${\cal P}_0$ which is not
residually finite. However, if we restrict to finitely generated groups
then there are only countably many of these in $\cal S$. In \cite{krvi}
it was shown that there are uncountably many finitely generated groups
in ${\cal P}_0$, so these two classes are far from equal though it
is still unknown whether $\cal P$ and ${\cal P}_0$ are the same.

\section{Actions on products of graphs}

We have looked throughout this paper at actions of (usually finitely
generated) groups $G$ on hyperbolic graphs $\Gamma$, with particular reference
to when $\Gamma$ is a quasitree and/or locally finite, or even of
bounded valence. The type of actions we have most been interested in are
those with unbounded orbits, or with finite stabilisers or metrically
proper actions. We were often able to find strong restrictions
on the structure of $G$, but in some cases it might be argued that the
conditions are too strong, for instance where groups containing $\Z^2$
were forbidden. A wider collection of metric spaces for a group to act
on by isometries are finite products of hyperbolic graphs and we investigate
this here in the last two sections. These spaces are well behaved in general
but are almost never hyperbolic.

We take a finite number of connected
graphs $X=\Gamma_1\times \ldots \times \Gamma_k$, each
equipped with their own path metric $d_i$, and turn $X$ into a metric space
using the $l_2$ metric. (We can also use the $l_1$ or $l_\infty$ metric
but we will see that the $l_2$ metric gives the richest structure.)
Hence for ${\bf x}=(x_1,\ldots ,x_k)$ and
${\bf y}=(y_1,\ldots ,y_k)\in X$ we have
\[d({\bf x},{\bf y})^2=\sum_{i=1}^kd_i(x_i,y_i)^2.\]
Consequently $X$ is a geodesic metric space because each factor is.
Note that $Isom(\Gamma_1)\times\ldots \times Isom(\Gamma_k)$ naturally embeds
in $Isom(X)$ by letting $(g_1,\ldots ,g_k)$ send $(x_1,\ldots ,x_k)$
to $(g_1(x_1),\ldots ,g_k(x_k))\in X$.
Suppose that we have a group $G$ acting on $X$ by isometries. We say
that $G$ {\it preserves factors} if the image of $G$ in $Isom(X)$
under the action lies in this subgroup. We can then think of $G$ as acting
separately by isometries on each factor $\Gamma_i$, so that we have a
homomorphism from $G$ to $Isom(\Gamma_1)\times \ldots \times Isom(\Gamma_k)$
and thus to each $Isom(\Gamma_i)$ separately by projection.
In particular any element
$g\in G$ can be written $g=(g_1,\dots ,g_k)$ for $g_i\in Isom(\Gamma_i)$.

However we can also use the graph structure on each $\Gamma_i$ to turn
$X$ into a $k$ dimensional complex, where a vertex of $X$ is
$(v_1,\ldots ,v_k)$ for $v_i$ a vertex of $\Gamma_i$, edges in $X$ join  
$(v_1,\ldots,v_i,\ldots,v_k)$ to $(v_1,\ldots,w_i,\ldots,v_k)$ for some $i$,
where $v_i-w_i$ is an edge in the graph $\Gamma_i$, and so on. We can then
consider the group $Aut(X)$ of combinatorial automorphisms, where this
structure is preserved. Any combinatorial automorphism is also an isometry
of $X$. Just as in the case for isometries we have that
$Aut(\Gamma_1)\times\ldots \times Aut(\Gamma_k)$ naturally embeds
in $Aut(X)$.
We also have that if a group $G$ acts on $X$ by combinatorial automorphisms
and it preserves factors then it acts on each factor by graph automorphisms.

If a group $G$ acts on $X$ by combinatorial automorphisms/isometries
then any subgroup $H$ clearly does too and in particular any finite index
subgroup. But alternatively,
if a group $H$ acts on $X$ by combinatorial automorphisms/isometries and $H$
is an index $i$ subgroup of $G$ then we can induce an action of $G$ on
$X^i$ by permuting factors, which will also be by combinatorial
automorphisms/isometries. Thus we obtain various properties invariant
under commensurability because $H$ has an action
by automorphisms on a finite product of graphs, each of which
comes from a particular class of graphs (eg
hyperbolic graphs, trees, locally finite graphs etc) if and only $G$
does. Moreover this  action of $H$ is metrically (or topologically)
proper, has unbounded
orbits etc if and only if the corresponding action of $G$ has this
property.

However this does not take into account actions that preserve factors.
In the case where $G$ acts by combinatorial automorphisms
of $X$ without preserving factors then it permutes some of these
graphs. Of course it can only permute graphs which are
isomorphic, but in any event we will always have a finite index subgroup
$H$ of $G$ which preserves factors and for the properties that concern
us, it is enough to work with $H$ instead of $G$.

Things can be somewhat different if $G$ is known
to act on $X$ by isometries but not necessarily automorphisms (at least
if the $l_2$ metric is used).
For instance if we take $k$ copies of the simplicial line $L$ then
$L^k$ with its natural cube complex structure has an automorphism group
which is virtually a copy of $\Z^k$, whereas $Isom(L^k)=Isom(\R^k)$ is
a much bigger group. Fortunately 
the wide ranging generalisations of de Rham theory
in \cite{foly} show that this is the only exception, as we obtain
the following:
\begin{prop} \label{ismprd}
Suppose that $X=\Gamma_1\times \ldots \times \Gamma_k$ is a finite product of
connected
graphs, where each $\Gamma_i$ has the induced path metric and $X$
has the $\ell_2$ product metric. Suppose that $G$ is any group
acting by isometries on $X$.\\
(i) If no graph $\Gamma_i$ is the simplicial line then $G$
has a finite index subgroup $H$ which preserves factors
(and thus $G$ acts by graph automorphisms on each factor $\Gamma_i$
which is not a circle).\\
(ii) If $\Gamma_{j+1},\ldots,\Gamma_k$ are those factors which are the
simplicial line $L$ then $G$ has a finite index subgroup $H$ which
preserves each factor in the decomposition
$\Gamma_1\times \ldots \times \Gamma_j\times \R^{k-j}$ where $H$ acts
on each $\Gamma_1,\ldots ,\Gamma_j$ by isometries (and therefore
by graph automorphisms on each factor that is not a circle) and by
Euclidean isometries on $\R^{k-j}$.
\end{prop}

Thus if we assume that all group actions are
by combinatorial automorphisms, Proposition \ref{ismprd}
shows us that the only cases acting isometrically and not covered by this
would have a de Rham factor $\R^k$ where (on dropping to a finite index
subgroup)
the action is by Euclidean isometries on this factor and by 
graph automorphisms on all the other factors
(or there would be a graph isomorphic to a circle
in the product, which is compact anyway). Nevertheless we will find that
allowing arbitrary isometric actions gives us a wider range of groups,
as discussed in the last section.

Thus when considering these various properties of actions, by dropping
down to a finite index subgroup we can always
assume for actions by automorphisms (and also for actions by isometries
if we include the de Rham factor) 
that if a group $G$ acts on $X$ in this way then it does so preserving
factors. But now we can ask: if $G$ does preserve factors
and the overall action of $G$ on $X$ has
a particular property then is this still true for the action
of $G$ on each factor? It is easily checked that $G$ acts with bounded orbits
on $X$ if and only if it acts with bounded orbits on each factor.
However, more interesting is the case where $G$ acts properly on $X$, as this
does not mean that $G$ acts properly on the factors even if it preserves them
(and even if there is no de Rham factor). Note that if our group action
on $X$ is by automorphisms then acting topologically properly is equivalent
to every vertex in $X$ having a finite stabiliser, just as in the case
for a single graph.

\section{Proper actions on products of quasitrees}

In this section we consider the question: do we have something akin to
a Tits alternative for groups acting properly on the finite product
of quasitrees or trees? We clearly cannot say anything in the most general
situation, by
recalling that any group acts topologically properly
and by automorphisms
on a (bounded) quasitree which is a graph (though not a locally finite one).
Thus we now have the
option of insisting that the action is metrically proper and/or restricting
to quasitrees which are locally finite graphs (in which case the two types of
proper are equal) and/or examining the case where the quasitrees are
simplicial trees. We can also wonder what happens if acting by automorphisms
is weakened to acting by isometries.

In the case where we insist that all factor graphs are simplicial trees,
the resulting product is a CAT(0) cube complex which will be finite
dimensional (but not locally finite unless the trees are). We have
versions of Tits alternatives for groups acting suitably on
these spaces in \cite{nibw} and \cite{caprr}
(see also \cite{ily} for a version in the case of median spaces).
We start by giving a variant in the case where we have products of trees. 

\begin{thm} \label{ttst}
  Let $G$ be any group acting topologically properly by automorphisms on a
  product of $k$ simplicial trees (which need not be locally finite).
  Then
either a finite index subgroup of $G$ has a action of general type on one
of these trees (so that $G$ contains $F_2$) or $G$ has a finite
index subgroup $H$ which is of the form $H/K\cong \Z^n$ for some
$0\leq n\leq k$,
where $K$ is a normal subgroup of $H$ which is locally finite.

In particular if $G$ contains no subgroup which is infinite but locally
finite, or even if $G$ contains no locally finite, infinite subgroup
having a normaliser of finite index, then either $F_2\leq G$ or $G$ is
virtually $\Z^n$ for some $n$.
\end{thm}

\begin{proof}
Let this action be on the product
$X=T_1\times \ldots \times T_k$ of trees.
We drop to a finite index subgroup of $G$, which we also call $G$,
acting topologically properly on $X$ and preserving factors, so that $G$ will
also act by automorphisms on each $T_i$.

We now have $k$ separate actions of $G$, each on a simplicial tree.
If any of these actions are of general type then we have the first case in the
statement of the theorem. Otherwise all remaining actions must be bounded,
parabolic, lineal or quasi-parabolic.
Renumber these trees so that the
lineal actions of $G$ are on $T_1,\ldots , T_l$ say, where we have an
invariant pair of ends. For each
$1\leq i\leq l$, 
we can take a subgroup $G_i$ of index 1 or 2 in $G$ where both ends are
fixed pointwise. This allows us to define a
Busemann function on $G_i$, which is a homogeneous
quasi-morphism of $G_i$ where the elements
acting elliptically on $T_i$ form the quasi-kernel. But for a tree this
function is actually a homomorphism from $G_i$ to $\Z$, which can be made onto
by rescaling. In each case the kernel
$K_i$ consists of the elements of $G_i$ 
acting elliptically on $T_i$ so we can now drop
to $H:=G_1\cap \ldots \cap G_l$ which also has finite index in $G$. Note
throughout that on dropping to a finite index
subgroup, we still have an action that preserves factors and moreover
it preserves the hyperbolic type of action on each factor. 

Now we consider the quasi-parabolic actions of $H$ (which are the same
as those of $G$) on the
trees $T_{l+1},\ldots ,T_{l+p}$ say. For each $l+1\leq i\leq l+p$, we have
an action fixing a point on the boundary, again giving us a
Busemann function on $H$ in the form of  a homomorphism from $H$ to
$\Z$ with the elements
acting elliptically on $T_i$ forming the kernel, as we have a tree.
We label these
$\theta_{l+1}, \ldots ,\theta_{l+p}$ in line
with the homomorphisms obtained from the lineal actions.

Thus we now have a homomorphism $\Theta=(\theta_1,\ldots ,\theta_{l+p}):H
\rightarrow \Z^{l+p}$ which may not
be onto but we can regard it as
a surjective homomorphism from $H$ to $\Z^n$ for $0\leq n\leq l+p$.
Let the kernel be called $K$: it consists of all elements of $H$
whose action
on each of the trees $T_1,\ldots ,T_{l+p}$ is elliptic. But all remaining
actions of $H$ on the other trees are parabolic or bounded, which means
that every element acts elliptically here and so every element of $K$ acts
elliptically on all of the factors.

Now take $S$ to be any finitely generated subgroup of $K$. Then
$S$ acts topologically properly on $X$ too, but for each tree $T_i$ we
have that every element of $S$ is acting elliptically. Thus this action
of $S$ on $T_i$ has bounded orbits and indeed by \cite{ser} a global fixed
point $v_i$ (which we can assume to be a vertex of $T_i$ by subdividing if
necessary). Thus $S$ fixes a vertex in each factor tree and thus
a vertex in $X$, thus this action can only be topologically proper
if $S$ is finite and therefore $K$ is locally finite.

Our subgroup $K$ of elements acting elliptically on every factor is
normal in $H$ and hence the normaliser of $K$ has finite index in
the original group $G$. Thus if $G$ has no infinite, locally finite
subgroups with a finite index normaliser then $K$ must be finite.
Now $H/K$ is isomorphic to $\Z^n$ and so $H$ and $G$ are finitely
generated, but a finitely generated group
which is finite-by-abelian is also abelian-by-finite, hence $H$ and $G$ are
virtually $\Z^n$.
\end{proof}

Note that we really do need the action to be by automorphisms here: in
\cite{bh} II.7.13 we have an example of $\Q$ acting properly on a 
2 dimensional cube complex which is locally finite, being the
product of a locally finite simplicial tree $T$ with $\R$. But $\Q$ is
definitely excluded
from the statement in Theorem \ref{ttst}. This ambiguity is resolved
by the fact that the action of $\Q$ on the $\R$ factor is by isometries
but not automorphisms.

Note also that if we ``Ripsify'' this action by replacing $\R$ (or even
$\Q$) with the Rips graph $\Gamma_1(\R)$ then the action of $\Q$ by
automorphisms on the product $T\times \Gamma_1(\R)$ of quasitrees is even
metrically proper. Thus the conclusion of Theorem \ref{ttst} cannot
hold if we replace ``simplicial trees'' in the hypothesis with
``graphs that are all quasitrees'', even if we strengthen acting topologically
properly to acting metrically properly. By now though it should be no surprise
that we can establish this conclusion if our quasitrees
are locally finite (whereupon acting metrically and topologically properly
are of course the same thing).

\begin{co} \label{ttsq}
Let $G$ be any group acting properly by automorphisms on a
product of $k$ graphs which are all locally finite quasitrees. Then
either a finite index subgroup of $G$ has a action of general type on one
of these quasitrees (so that $G$ contains $F_2$) or $G$ has a finite
index subgroup $H$ which is of the form $H/K\cong \Z^n$ for some
$0\leq n\leq k$,
where $K$ is a normal subgroup of $H$ which is locally finite.

In particular if $G$ contains no subgroup which is infinite but locally
finite, or even if $G$ contains no locally finite, infinite subgroup
having a normaliser of finite index, then either $F_2\leq G$ or $G$ is
virtually $\Z^n$ for some $n$.
\end{co}

\begin{proof}
Let $X$ now be $Q_1\times\ldots \times Q_k$. 
We still have that any element acting on one of these
quasitrees $Q_i$  is either elliptic
or loxodromic and again the action of $G$ on each $Q_i$ 
is either bounded, parabolic, lineal or quasi-parabolic. We can now
follow the proof of Theorem \ref{ttst}, using
Theorem \ref{two} for the quasitrees where the action is lineal
as they are locally finite, regardless of whether $G$ is finitely generated.
For each of the quasi-parabolic
actions we mentioned that all orbits are quasi-convex. Thus
again we apply Proposition
\ref{os} so that the constructions with quasitrees in Section 6 work here
too. This gives us the same conclusion as for Corollary \ref{bndd},
namely that $H$ acts coboundedly, by graph automorphisms and with unbounded
orbits
on a bounded valence quasitree which is quasi-isometric to a
bounded valence tree $T$ which is bushy.
We can now apply \cite{msw} in the form of Theorem \ref{summ} to get
an action of $H$ on some bounded valence tree $T'$ where there is a
coarse $H$-equivariant quasi-isometric embedding $q_i$
from the action of $H$
on $T'$ to the original action of $G$ on $Q_i$

It is now easily checked that, on replacing $X=Q_1\times \ldots \times Q_i\times
\ldots \times Q_k$ with
$Y=Q_1\times \ldots \times T'\times
\ldots \times Q_k$ then the map
$(id,\ldots ,q_i,\ldots ,id)$ from $X$ to $Y$
is also a coarse $H$-equivariant quasi-isometric
embedding. In particular the property of $H$ acting properly on $X$ passes
over to the action of $H$ on $Y$. We now do this in turn for all of the
quasitrees where the action of $H$ is either lineal or quasi-parabolic.

Thus we obtain our homomorphism from $H$ to $\Z^{l+p}$ as before with
the kernel $K$ acting elliptically on each factor. Again $K$ is
locally finite, by replacing \cite{ser} with
Corollary \ref{serc}.
\end{proof}

We also have:
\begin{co} \label{smgp}
  If $G$ is a finitely generated group acting properly by
  automorphisms (respectively isometries)
  on a finite product of graphs which are all
  locally finite quasitrees then $G$ acts properly by automorphisms
(respectively isometries)
  on a finite product of bounded valence trees.
\end{co}
\begin{proof} We can assume that $G$ acts by preserving factors because
  we can drop down to a finite index subgroup that does this and then
  back up to $G$ at the end. In the case where $G$ acts by automorphisms,
  we have an action of $G$ on each of the factor
  graphs $Q_1,\ldots ,Q_k$ say, which is also by automorphisms. As before
  we can assume that none of these actions
  have bounded orbits, because if so we could remove this graph
  from the product and $G$ would still act properly on the new product.
  Thus we can now apply Theorem \ref{summ}
  to each of these actions (regardless of their hyperbolic type, as $G$
  is finitely generated)
  to obtain an action of $G$ on a bounded valence tree $T_i'$ from the
  action of $G$ on $Q_i$. Putting these together gives us an action by
  automorphisms of $G$ on the product $T'_1\times \ldots \times T'_k$
  of bounded valence trees, and this will still be a proper action as
  in Corollary \ref{ttsq} above.

  For an action of $G$ by isometries, we can drop down to a finite
  index subgroup and say that without loss of generality $G$ acts
  on the factor graphs by automorphisms (or a factor graph is a circle
  but we can remove this without altering properness of the action)
  and on the de Rham factor by isometries. But from the above we can replace
  the factor graphs with actions of trees, so we now have $G$ acting
  properly by isometries on the product of locally finite trees and $\R^m$,
  so by isometries on the product of locally finite trees.
\end{proof}

We might ask about the existence of a finitely generated group acting
properly by isometries on a finite product of graphs which are all
locally finite quasitrees but which has no such action by automorphisms.
We will now examine a much stronger counterexample. To do this, we
first look at Bestvina, Bromberg and Fujiwara's property (QT).
A finitely
generated group is said to have (QT) if it acts by isometries
on a finite product of
quasitrees such that the orbit map is a quasi-isometric embedding. 
(In this definition a quasitree is always a graph equipped with its
path metric which is quasi-isometric to a simplicial tree but which need
not be locally finite.) This is a strong definition and fits
into our approach here because we mentioned in Section 2 that
that if a finitely generated group $G$
acts on a metric space
$X$ by isometries then the orbit map (with respect to an arbitrary basepoint
$x_0\in X$) is a quasi-isometric embedding if and only if there exists
some coarse $G$-equivariant quasi-isometric embedding from $G$
(with its word metric) to $X$. We also mentioned that 
it implies the action is metrically proper.

In \cite{bebf} it is shown that mapping
class groups and all residually finite hyperbolic groups have (QT),
which is not true in general in either case if the quasitrees
are locally finite. For the second case we can
take a hyperbolic group with property (T) which will have no proper
action by isometries on a finite product of locally finite quasitrees
by Corollary \ref{smgp} (indeed any such action will have finite orbits).
In the case of (most) mapping class groups, we have that
they possess no proper action by isometries on a product of locally
finite quasitrees. However, as for unbounded actions, we remark that
even in the case of a product of locally finite simplicial trees this
is unknown (apart from the case of genus 2 where the answer is known to
be yes). Indeed, it is a well known open question whether some finite index
subgroup splits as an HNN extension or amalgamated free product. 

It is also a consequence of \cite{drnjan} that every Coxeter group has
(QT). Note that (QT) is preserved by finite index supergroups and descends
to undistorted finitely generated subgroups. Other results on groups,
including 3-manifold groups, having property (QT) can be found
in \cite{bei1} and \cite{bei2}.

Our group of interest is the recent group of
Leary - Minasyan in \cite{lrmn} given by the finite presentation
\[G=\langle t,a,b\,|\,[a,b],ta^2b^{-1}t^{-1}=a^2b,tab^2t^{-1}=a^{-1}b^2\rangle. \]
This is a CAT(0) group and is an HNN extension of $\Z^2=\langle a,b\rangle$
with edge subgroups having index 5. In particular it acts properly and
cocompactly by isometries on
the product of the regular simplicial tree $T_{10}$ and $\R^2$.
Thus $G$ acts by isometries
on a finite product of locally finite trees and the orbit map is a
quasi-isometry. Thus it would certainly appear to have a very strong form
of the property (QT). Nevertheless we will now show that $G$ does not have
(QT). This apparent contradiction is because the product of quasitrees
in the definition of (QT) is equipped with the $l_1$ product metric.

\begin{thm} \label{ashl}
  The Leary - Minasyan finitely presented group
  \[G=\langle t,a,b\,|\,[a,b],ta^2b^{-1}t^{-1}=a^2b,tab^2t^{-1}=a^{-1}b^2\rangle \]
  acts geometrically by isometries on a finite product of locally finite
  simplicial trees (using the path metric), 
   with the $l_2$ product metric. However it
  does not possess any metrically proper action by isometries
  on any finite product of graphs (using the path metric)
  which are quasitrees, with the $l_1$ or $l_\infty$ metrics. In particular
  $G$ does not have property (QT). Nor
  does $G$ possess any action by automorphisms with finite stabilisers
  on any finite product of locally finite graphs (using the path metric)
  which are quasitrees.
\end{thm}
\begin{proof}
  The first statement is shown above because the product
  $T_{10}\times\R\times\R$ of locally finite trees is $T_{10}\times \R^2$
  provided we use the $l_2$ product metric.

  For the second statement we will show below in Theorem \ref{ashla}
  that, on setting
  $A=\langle a,b\rangle\leq G$, for any finite
    index subgroup $H$ of $G$ and any action of $H$ by isometries
    on any hyperbolic space, no element of the subgroup
    $H\cap\langle a,b\rangle\cong \Z^2$ can act loxodromically.
    Assuming this here, we now consider the equivalent version of 
    Proposition \ref{ismprd} for the $l_1$ and $l_\infty$ product metrics
    to show that a de Rham factor does not occur.

\begin{prop}   \label{lonpr} 
Suppose that $X=\Gamma_1\times \ldots \times \Gamma_m$ is a finite product of
connected
graphs, where each $\Gamma_i$ has the induced path metric and $X$
has the $l_1$ or the $l_\infty$
product metric. Suppose that $G$ is any group
acting by isometries on $X$. Then $G$
has a finite index subgroup $H$ which preserves factors
(and thus $H$ acts by graph automorphisms on each factor $\Gamma_i$
which is not a circle or the simplicial line).
\end{prop}
\begin{proof}
  On giving a graph the induced path metric, it becomes not just a
metric space which is geodesic
but a path connected metric space which is locally uniquely geodesic.
Thus if $X$ is given the $l_\infty$ metric then this result is a direct
consequence of \cite{malw}, which is established by using the fact that for
certain directions in the product space, geodesics are unique if they are
unique in the factors. This is certainly true in an $l_1$ product space
if we travel ``horizontally/vertically'', thus
we can mimic his proof for this (easier) case:
\begin{lem} \label{unig}
  Let $X=X_1\times\ldots\times X_m$ have the $l_1$ product
  metric $d_X$
  for geodesic metric spaces $(X_1,d_1),\ldots ,(X_m,d_m)$.
Suppose we have  two points ${\bf x},{\bf y} \in X$ of the form
\[{\bf x}=(x_1,\ldots ,x_{i-1},x_i,x_{i+1},\ldots ,x_m)\mbox{ and }{\bf y}=
  (x_1,\ldots ,x_{i-1},y_i,x_{i+1},\ldots ,x_m)\]
where $1\leq i\leq m$. Then any geodesic between $\bf x$ and $\bf y$
only varies in the $i$th coordinate. Conversely if we have two
points ${\bf x},{\bf y} \in X$ which differ in at least two coordinates
then there is more than one geodesic from $\bf x$ to $\bf y$ in $X$.
\end{lem}
\begin{proof}
  Suppose we have a (unit speed) geodesic
  $\boldsymbol{\gamma}:[0,d]\rightarrow X$ from $\bf x$ to $\bf y$ where
  $d=d_X({\bf x},{\bf y})$. If there is $j_0\neq i$ and $t\in (0,d)$
  such that the $j_0th$ coordinate of
  $\boldsymbol{\gamma}(t)=(t_1,\ldots ,t_m)$ is not equal to $x_{j_0}$ then
  \begin{eqnarray*} d_i(x_i,y_i)=
    d_X({\bf x},{\bf y})&=&d({\bf x},\boldsymbol{\gamma}(t))
      +d(\boldsymbol{\gamma}(t),{\bf y})\\
                        &=&\Sigma_{j=1,j\neq i}^m \left(d_j(x_j,t_j)
                            +d_j(t_j,x_j)\right)
                            +d_i(x_i,t_i)+d_i(t_i,y_i)\\
                        &\geq& d_{j_0}(x_{j_0},t_{j_0})+d_{j_0}(t_{j_0},x_{j_0})+
                               d_i(x_i,t_i)+d_i(t_i,y_i)\\
                          &>& d_i(x_i,y_i)\mbox { as } x_{j_0}\neq t_{j_0}
\end{eqnarray*}
which is a contradiction.

Now suppose we have two points ${\bf x}=(x_1,\ldots ,x_m)$ and
${\bf y}=(y_1,\ldots ,y_m)$ in $X$ which differ in (without loss of generality)
at least the first two coordinates. Take unit speed
geodesics $\gamma_i$ in each $X_i$ running from $x_i$ to $y_i$ and set
\[\boldsymbol{\delta}_1(t)=(x_1,\gamma_2(t),\ldots ,\gamma_m(t)),
  \boldsymbol{\delta}_2(t)=(\gamma_1(t),y_2,\ldots ,y_m))\]
and
\[\boldsymbol{\delta}_3(t)=(\gamma_1(t),x_2,\gamma_3(t),\ldots ,\gamma_m(t)),
  \boldsymbol{\delta}_4(t)=(y_1,\gamma_2(t),y_3,\ldots ,y_m)).\]
Then it is easily checked that following $\boldsymbol{\delta}_1$ then
$\boldsymbol{\delta}_2$
and also following $\boldsymbol{\delta}_3$ then $\boldsymbol{\delta}_4$
are both geodesics
from $\bf x$ to $\bf y$ and they are distinct as $x_1\neq y_1$
and $x_2\neq y_2$.
\end{proof}

Now we return to the setting in Proposition \ref{lonpr}. Each $\Gamma_i$
is a geodesic metric space (which we
assume without loss of generality is not a single point). Thus given
any isometry $g$ of $X$, take any
point ${\bf x}=(x_1,\ldots ,x_i,\ldots ,x_m)$ in $X$ and let $z_i$ be another
point in $\Gamma_i$ near $x_i$ such that there is only one geodesic
$\gamma_i$ in $\Gamma_i$ from $x_i$ to $z_i$. Then by Lemma \ref{unig}
we have that 
$\boldsymbol{\gamma}(t)$, which is equal to
$(x_1,\ldots ,\gamma_i(t),x_{i+1},\ldots ,x_m)$,
is the unique geodesic in $X$ between $\bf x$ and
${\bf z}=
(x_1,\ldots ,x_{i-1},z_i,x_{i+1},\ldots ,x_m)$ and so under $g$
it must map to
a unique geodesic between $g({\bf x})$ and $g({\bf z})$, by considering
$g^{-1}$. Thus by Lemma \ref{unig} again
these two points (which are not the same because $g$ is
a bijection) differ in only one coordinate, say the $j$th.

Now take $y_i$ to be an arbitrary point in $\Gamma_i$. Given any geodesic
in $\Gamma_i$ from $x_i$ to $y_i$, we can split it into a finite number
of subgeodesics, with some overlap that is more than a point, where each
subgeodesic is the unique geodesic between its own endpoints. Thus the image
of each of these subgeodesics under $g$ is a subset of $X$ where only
one coordinate varies, and as these subgeodesics overlap this will always
be the $j$th coordinate. In other words given $1\leq i\leq m$ we have
$1\leq j\leq m$ such that for fixed
$x_1\in \Gamma_1,\ldots ,x_{i-1}\in \Gamma_{i-1},x_{i+1}\in\Gamma_{i+1},\ldots ,
x_m\in\Gamma_m$ there is
$y_1\in \Gamma_1,\ldots ,y_{j-1}\in \Gamma_{j-1},y_{j+1}\in\Gamma_{j+1},\ldots ,
y_m\in \Gamma_m$ such that we have a function $f:\Gamma_i\rightarrow \Gamma_j$
with

\[g(x_1,\ldots ,x_{i-1},x,x_{i+1},\ldots ,x_m)=
  (y_1,\ldots ,y_{j-1},f(x),y_{j+1},\ldots ,y_m)\]
for all $x\in \Gamma_i$. Moreover $f$ is a
bijection and an isometry because $g$ is. By letting $g$ vary over
all the isometries in $G$ we obtain our result.
\end{proof}                              
We can now finish the proof of Theorem \ref{ashl}. If $G$ is acting
by isometries and metrically properly on a finite product of graphs
which are quasitrees with the $l_1$ or $l_\infty$ product metric then
we now know we have a finite index subgroup $H$ (also acting metrically
properly) which preserves factors
and so acts by isometries on each factor. Therefore this is also true
when we restrict the action to $H\cap A$. But as isometries of quasitrees
can only be loxodromic or elliptic, the action of $H\cap A$ is purely
elliptic on each quasitree and hence is bounded, by Corollary \ref{serc}
because $H\cap A$ is finitely generated (or just directly as
$H\cap A\cong\Z^2$).
Therefore its action on $X$
is also bounded, so we have an infinite subgroup of $G$ acting on $X$
with bounded orbits, meaning that the action of $G$ is not
metrically proper.
\end{proof}

In order to show that certain elements of $G$ cannot act as loxodromic
elements under any action on a hyperbolic space, we consider stable
translation length. We have already seen that we can have actions of
$\Z^2=\langle a,b\rangle$ on a hyperbolic space where every non identity
element is loxodromic, such as $a\mapsto a+1$ and $b\mapsto b+\sqrt{2}$
on $\R$, as well as related actions of this form. We now show that this
is the only type of action with every non identity element loxodromic,
at least at the level of
stable translation length $\tau$. This was defined in Section 3: recall
that for an isometric
action of an arbitrary group $G$ on an arbitrary metric space
$(X,d)$, it satisfies $\tau(g)\leq d(g(x),x)$ as well as:
\begin{eqnarray*}\bullet &&\tau(g)=0\mbox{ if and only if $g$ is not
                            loxodromic} \\
\bullet &&\tau(g^n)=|n|\tau(g)\mbox{ for all }n\in\Z\\
  \bullet &&\tau(gh)\leq \tau(g)+\tau(h)\mbox{ {\bf if} $G$ is abelian}.
\end{eqnarray*}
In particular, if we set $G$ equal to the free abelian group $\Z^m$ we have
that the function $\tau:G\rightarrow [0,\infty)$ is a $\Z$-seminorm and
further is a $\Z$-norm if and only if all non identity elements act
loxodromically. We also have the concepts of an $\R$-norm and an $\R$-seminorm,
defined similarly for a vector space over $\R$ using the scalars. The crucial
point here for actions of $\Z^2$ is: 
\begin{thm} \label{stz}
  If $\Z^2$ acts by isometries on an arbitrary hyperbolic space $X$ then
  there is a homomorphism $\theta:\Z^2\rightarrow \R$ such that
  \[| \theta (g)|=\tau(g)\mbox{ for all }g\in\Z^2.\]
\end{thm}
\begin{proof} First suppose that all elements act loxodromically so that
  $\tau$ is  a $\Z$-norm on $\Z^2$, where we write the element $a^mb^n$
  additively as $(m,n)$ also.
  Then $\tau$ naturally extends to $\Q^2$ and to $\R^2$
  by taking limits, which gives us an $\R$-seminorm. Suppose that this is
  actually an $\R$-norm, that is $v\in\R^2\setminus \{(0,0)\}$ implies
  $\tau(v)\neq 0$. As all $\R$-norms on $\R^2$ are equivalent, this gives
  us $c>0$ such that for all $(x,y)\in \R^2$ we have
  $c(|x|+|y|)\leq \tau(x,y)$ and so this certainly holds
  in $\Z^2$ as well. In other
  words, for any point $x_0\in X$ and element $a^mb^n$ we have
  \[c(|m|+|n|)\leq\tau(m,n)\leq d(a^mb^n(x_0),x_0).\]
  But we also have
  \[d(a^mb^n(x_0),x_0)\leq K(|m|+|n|)\mbox{ for }
    K=\mbox{max}\{d(a(x_0),x_0),d(b(x_0),x_0)\}.\]
  This means that $(m,n)\mapsto a^mb^n(x_0)$ is a quasi-isometric
  embedding of $\Z^2$ in the hyperbolic space $X$ which is a contradiction.

  Thus $\tau$ is in fact an $\R$-seminorm and there is a one dimensional
  vector subspace $S$ of $\R^2$ on which $\tau$ is zero. Thus $\tau$ factors
  through to an $\R$-norm on the quotient one dimensional vector space
  $\R^2/S$ but the only $\R$-norms on a one dimensional real vector
  space $V$ are $v\mapsto |\mu\alpha(v)|$ for some $\mu\in \R\setminus \{0\}$,
  a non zero linear map $\alpha:V\rightarrow \R$ and $|\cdot |$ the usual
  modulus on $\R$.

  Thus if $\theta$ is the composition of the linear maps
  $r\mapsto r+S$ where $r\in \R^2$ and $v\mapsto \alpha(v)$ where
  $v\in \R^2/S$, we have for appropriate $\mu$ and $\alpha$ that $\theta$
  is linear and $\tau(r)=|\theta (r)|$.
  
  If there is an element in $\Z^2$ that is not loxodromic, thus is zero
  under $\tau$, then we can still define $\tau$ on $\R^2$ to get an
  $\R$-seminorm but it is not an $\R$-norm on $\R^2$. In this case either
  no elements are loxodromic and $\theta\equiv 0$ or we can quotient out
  by a one dimensional subspace and argue just as before.
\end{proof}

We can now apply this to the Leary - Minasyan group $G$ with subgroup
$A=\langle a,b\rangle\cong \Z^2$.
\begin{thm} \label{ashla}
  Given any finite index subgroup $H$ of $G$ and any isometric action of $H$
  on any hyperbolic space $X$, no element of $A\cap H\cong\Z^2$
  can act as a loxodromic isometry.
\end{thm}
\begin{proof}
  We have $ta^5t^{-1}=a^3b^4$ and $tb^5t^{-1}=a^{-4}b^3$ holding in $G$. This
  implies by induction on $n$ that
  \[t^na^{5^n}t^{-n}=a^{\alpha_n}b^{\beta_n}
    \mbox{ and }t^nb^{5^n}t^{-n}=a^{\gamma_n}b^{\delta_n}\]
  where
  \[\sma{cc} \alpha_n&\beta_n\\\gamma_n&\delta_n \fma
    =\sma{rc} 3&4\\-4&3\fma^n.\]
  Take $K\in \N$ such that $s:=t^K$, $c:=a^K$ and $d:=b^K$ are all
  in $H$. Then we have
  \begin{eqnarray*}
    sc^{5^K}s^{-1}&=&c^{\alpha_K}d^{\beta_K}\\
    \mbox{ and }\qquad sd^{5^K}s^{-1}&=&c^{\gamma_K}d^{\delta_K}
  \end{eqnarray*}
  holding in $H$. We now consider what this means for the stable
  translation length function $\tau$ on $H$. As $\tau$ is conjugation
  invariant, on applying Theorem \ref{stz} to $H\cap A$ we get a
  homomorphism $\theta:H\cap A\rightarrow \R$ which, on setting
  $x:=\theta(c)$ and $y:=\theta(d)$, satisfies the two equations
  \begin{eqnarray*}
    \pm 5^Kx&=&\alpha_Kx+\beta_Ky\\
    \pm 5^Ky&=&\gamma_Kx+\delta_Ky.
  \end{eqnarray*}

  Thus by rescaling the hyperbolic metric on $X$ and hence $x$ and $y$,
  we can assume that they are both rational and so the homomorphism
  $\theta$ has discrete image in $\R$. This means that there is some
  non identity element $c^md^n\in H\cap A$ such that $\theta(c^md^n)=0$.
  But note that on $H\cap A$ we have $\theta$ is zero exactly when $\tau$
  is zero. Thus we also have $sc^md^ns^{-1}$ and indeed $sc^{m5^K}d^{n5^K}s^{-1}$
  mapping to zero under $\theta$. But the latter element is equal to
  \[c^{m\alpha_K+n\gamma_K}d^{m\beta_K+n\delta_K}\]
  and so $\theta$ and hence $\tau$ will be identically zero (meaning that
  no element in $H\cap A$ is loxodromic) unless
  \[(m,n) \mbox{ and }(m\alpha_K+n\gamma_K,m\beta_K+n\delta_K)\]
  are linearly related. But this is the same as saying that
  $\sma{c} m\\n\fma$ is an eigenvector of the matrix
  $\sma{cr} 3&-4\\4&3 \fma^K$, which would imply that this matrix has
  a rational eigenvalue. But this is the matrix of an infinite order
  rotation, so its eigenvalues are not even real.
\end{proof}

Having now seen restrictions on finitely generated
groups that can act properly by
automorphisms on a finite product of locally finite quasitrees, we
might ask: which (finitely generated)
groups are definitely known to have such actions? 
In the case of trees we have direct products of free groups and
Burger - Mozes groups amongst others. However if we
now consider which word hyperbolic groups have such an action, we
find that here things are wide open. Of course a virtually free group
has such an action by taking a single tree and we have shown in
Corollary \ref{smgp} that turning each tree into a quasitree allows no more
groups.

For products of trees/quasitrees, we now show that all of these possible
variants are the same in the case of hyperbolic groups.
\begin{thm} \label{hypg}

  Let $G$ be a word hyperbolic group and consider the statements:\\
  \hfill\\
  $G$ has a proper action by automorphisms/isometries on a finite
  product of locally finite/bounded valence/$d$-regular (for some $d$)
  graphs which are\\
  trees/quasitrees, equipped with the $l_1$/$l_2$/$l_\infty$ product metric.

  Then these statements are all equivalent, for any possible choice of
  the different options.
\end{thm}
\begin{proof}
  If we first take automorphisms, trees and the $l_2$ product metric
  in the above then the equivalence of these for locally finite, bounded
  valence and $d$ regular (for some $d$) trees  was shown in \cite{me1}
  Lemma 6.1.

  If $G$ is acting properly by automorphisms on a finite product
  of graphs which are locally finite quasitrees, equipped
  with the $l_2$ product metric, then Corollary \ref{smgp} tells us
  that we can assume our graphs are bounded valence trees (and consequently
  we can take them to be $d$-regular for some $d$ by the above). This works
  because $G$ acting by automorphisms
  on the product space means that we can drop to a finite index subgroup $H$
  (and induce the action back to $G$ at the end) which
  acts by automorphisms on each factor graph. Then $H$ will of course
  act by isometries on each factor graph and hence by isometries on the
  product, regardless of whether the $l_1$, $l_2$ or $l_\infty$ product
  metric is chosen. Moreover whether the action is proper will be
  independent of this choice too, thus Corollary \ref{smgp} also applies
  to the other product metrics when the action is by automorphisms.

  Now let us suppose that our action is by isometries with
  the $l_2$ product metric. Then Corollary \ref{smgp} also allows us to
  replace our product of locally finite quasitrees with an isometric action
  on a product of bounded valence trees.

  This also works for the $l_1$ and $l_\infty$ product metrics, whereupon
  we have an isometric action by $G$ (or a finite index subgroup, which
  we will also call $G$) on
  $T_1\times\ldots\times T_m\times\R^k$ which preserves factors and
  which acts by
  automorphisms on the locally finite trees $T_1,\ldots ,T_m$ (this includes
  the cases of  the $l_1$ or $l_\infty$ product metric where we can further
  say without loss of generality that the isometric action on $\R^k$ preserves 
  each component).

  If the action of $G$ is also proper on restricting to
  $X=T_1\times\ldots\times T_m$ then we are done. If not then there must
  exist a vertex ${\bf v}$ of $X$ with infinite stabiliser $S\leq G$. The
  action of $S$ on $X\times \R^k$ must also be proper, so $S$ acts properly
  on $\R^k$ (or on each factor $\R$ in the $l_1$ and $l_\infty$ cases)
  by Euclidean isometries and so must be virtually abelian, say with normal
  abelian subgroup $A$ of index $i$ in $S$. As $G$ is (non elementary)
  word hyperbolic and $S$ is an infinite subgroup, we can take an infinite
  order element $g\in S$ and $h\in G$ such that no positive powers of
  $g$ and $hgh^{-1}$ commute. Now $hgh^{-1}$ fixes the vertex $h({\bf v})\in X$
  but as the (1-skeleton of) $X$ is a locally finite graph, some positive
  power of $hgh^{-1}$ fixes the neighbouring vertices of $h({\bf v})$ and
  another power fixes all vertices distance two away as well, and so on.
  Thus there is an integer $K>0$
  with $hg^Kh^{-1}$ in $S$. But then $hg^{iK}h^{-1}$ and
  $g^i$ are both in $A$ and so commute, which is a contradiction.

  Thus in this case we do not need the $\R^k$ factor, so (taking an
  induced action if necessary) $G$ acts properly by automorphisms
  on a finite product of bounded valence trees. We have already noted that
  this can be taken to be $d$-regular (for some $d\geq 3$, though the action
  might well not be minimal on these regular trees) and that the action
  will be proper regardless of which of our three product metrics are used.
\end{proof}
{\it Note}: We used much less than the full force of $G$ being word
hyperbolic in this proof; namely finite generation and that property
of non commuting elements. Nevertheless it is instructive to see how
this proof breaks down in the case of the Leary - Minasyan group.
   
It is unknown whether any word hyperbolic group has such an
action other than the virtually free groups. This is a well known
question in the case of the closed hyperbolic surface group $S_g$
(for genus $g\geq 2$, or equal to 2 without loss of generality).
See \cite{me1} for some positive evidence in this direction.

We can thus ask the two following open questions, which look quite
different but nevertheless are equivalent by Theorem \ref{hypg}.\\
\hfill\\
{\bf Question 1}: Is there $d\geq 3$ such that the surface group $S_g$
has an action by automorphisms on a finite product of $d$-regular trees
with finite stabilisers?\\
\hfill\\
{\bf Question 2}: Can the surface group $S_g$ have an action by isometries
on some finite product of graphs which are locally finite quasitrees,
with the $l_2$ product metric?

We can also ask the same two open questions (again equivalent to each other)
for any word hyperbolic group in place of $S_g$, other than the virtually
free groups.

We think it would be strange, to say the least, if the list of finitely
generated groups shown to have (QT) grew ever longer whilst in the case
of a finite product of locally finite quasitrees, the question remains of
whether there is any word hyperbolic group, other than the virtually free
groups, having a proper isometric action on such a space.

\end{document}